\documentclass[a4paper,11pt]{article}
\title{Hydrodynamic limit of $N$-branching Markov processes}
\author{Jean B\'erard and Brieuc Fr\'enais}
\date{}

\usepackage[utf8]{inputenc}
\usepackage[T1]{fontenc}
\usepackage{fullpage}

\usepackage{enumitem}

\usepackage{amsmath}			
\usepackage{amssymb}			
\usepackage{mathrsfs}
\usepackage{graphicx}
\usepackage{subcaption}

\usepackage{stmaryrd}			

\usepackage{imakeidx}

\usepackage{hyperref}	

\usepackage{xcolor}

\newtheorem{prop}{Proposition}
\newtheorem{theo}{Theorem}
\newtheorem{lemma}{Lemma}
\newtheorem{coro}{Corollary}

\newtheorem{rem}{Remark}

\newenvironment{demo}[1][Proof. ]
{\vspace{-0cm}\noindent\textit{#1}}
{\hfill $\square$ \vspace{0.3cm}}

\renewcommand{\a}{\alpha}

\newcommand{\BB}{\mathcal{B}}
\renewcommand{\b}{\beta}

\newcommand{\CC}{\mathcal{C}}

\newcommand{\DD}{\mathcal{D}}

\newcommand{\dd}{\delta}

\newcommand{\E}{\mathbb{E}}
\newcommand{\EE}{\mathcal{E}}
\newcommand{\ee}{\varepsilon}

\newcommand{\FF}{\mathcal{F}}

\newcommand{\G}{\mathcal{G}}
\newcommand{\Gg}{\Gamma}
\newcommand{\g}{\gamma}

\renewcommand{\k}{\kappa}

\renewcommand{\L}{\mathcal{L}}
\newcommand{\Ll}{\Lambda}

\newcommand{\NN}{\mathcal{N}}

\renewcommand{\O}{\Omega}

\renewcommand{\P}{\mathbb{P}}

\newcommand{\Q}{\mathbb{Q}}

\newcommand{\R}{\mathbb{R}}

\newcommand{\s}{\sigma}

\renewcommand{\t}{\theta}

\newcommand{\V}{\mathbb{V}}

\newcommand{\X}{\mathcal{X}}

\newcommand{\Z}{\mathbb{Z}}

\newcommand{\z}{\zeta}

\newcommand{\1}{\mathbf{1}}

\renewcommand{\d}{\partial}

\DeclareMathOperator{\e}{e}

\newcommand{\dint}{\displaystyle\int}

\newcommand{\liml}{\lim\limits}

\newcommand\lsto{\preccurlyeq}

\newcommand{\supl}{\sup\limits}

\newcommand{\dsum}{\displaystyle\sum\limits}

\newcommand{\cv}{\longrightarrow}

\newcommand{\ent}[2]{\llbracket #1,#2\rrbracket}

\newcommand{\n}[2]{\left\|#1\right\|_{#2}}

\renewcommand{\Tilde}{\widetilde}

\newcommand{\thref}[1]{\hyperref[#1]{Theorem \ref{#1}}}
\newcommand{\defiref}[1]{\hyperref[#1]{Definition \ref{#1}}}
\newcommand{\figref}[1]{\hyperref[#1]{Figure \ref{#1}}}
\newcommand{\propref}[1]{\hyperref[#1]{Proposition \ref{#1}}}
\renewcommand{\eqref}[1]{\hyperref[#1]{(\ref{#1})}}
\newcommand{\chapref}[1]{\hyperref[#1]{Chapter \ref{#1}}}
\newcommand{\lemref}[1]{\hyperref[#1]{Lemma \ref{#1}}}
\newcommand{\cororef}[1]{\hyperref[#1]{Corollary \ref{#1}}}
\newcommand{\secref}[1]{\hyperref[#1]{Section \ref{#1}}}

\newcommand{\ie}{\textit{\text{i.e.}}}

\begin{document}

\maketitle

\abstract{We consider the behaviour of branching-selection particle systems in the large population limit. The dynamics of these systems is the combination of the following three components: (a) Motion: particles move on the real line according to a continuous-time Markov process; (b) Branching: at rate $1$, each particle gives birth to a new particle at its current location; (c) Selection: to keep the total number of particles constant, each branching event causes the particle currently located at the lowest position in the system to be removed instantly. 
Starting with $N \geq 1$ particles whose positions at time $t=0$ form an i.i.d. sample with distribution $\mu_0$, we investigate the behaviour of the system at a further time $t>0$, in the limit $N \to +\infty$.
Our first main result is that, under suitable (but rather mild) regularity assumptions on the underlying Markov process, the empirical distribution of the population of particles at time $t$ converges to a deterministic limit, characterized as the distribution of the Markov process at time $t$ conditional upon not crossing a certain (deterministic) moving boundary up to time $t$. Our second result is that, under additional regularity assumptions, the lowest  particle position at time $t$ converges to the moving boundary. 

These results extend and refine previous works done by other authors, that dealt mainly with the case where particles move according to a Brownian motion. For instance, our results hold for a wide class of Lévy processes and diffusion processes. Moreover, we obtain improved non-asymptotic bounds on the convergence speed.

}

\tableofcontents

\section{Introduction}

We consider a branching-selection stochastic process that consists in a finite family of $N \geq 1$ particles moving on the real line in continuous time. At rate $1$, each particle gives birth to a new particle at its current location (branching), and, in order to preserve a constant number of particles, the particle currently located at the lowest position\footnote{If several particles are simultaneously located at the lowest position, ties are broken arbitrarily.} is removed instantly (selection). Between such branching-selection events, particles move independently according to a continuous-time Markov process $X$ on $\R$.  

This process belongs to a family of models introduced by Brunet and Derrida \cite{brunet-derrida}, which then gave rise to a series of related works (see \secref{subsec:bib} for a bibliographical discussion). Among these models, the so-called $N$-branching Brownian motion ($N$-BBM), where the underlying Markov process is the standard Brownian motion on $\R$, fits into our setting. In \cite{nbbm,nonlocal-nbbm}, De Masi et al. studied the hydrodynamic limit of the $N$-BBM (fixed time $t$, with $N\to+\infty$) and showed that it can be described as the solution of a free boundary problem. Our main goal is to extend this result to the case of a general Markov process on $\R$. Accordingly, we call this process the $N$-branching Markov process ($N$-BMP). 

\subsection{Statement of the main results}

Let $X=(X_t)_{t\geq 0}$ be a time-homogeneous Markov process on $\R$ (in the sense of \cite{revuz-yor}, Chapter III) with transition kernel $(p_t)_{t\geq 0}$\index{$(p_t)_{t\geq 0}$: transition kernel of $X$}, defined on a filtered measurable space $(\O,\FF,(\FF_t)_{t\geq 0})$ equipped with a family of probability measures $(\Q_x)_{x\in\R}$. As usual, we denote by $\Q_{\mu_0}$ the probability measure on $(\O,\FF)$ in the case where the distribution of $X_0$ is a certain probability measure $\mu_0$\index{$\mu_0$: initial distribution of the particles}. We denote by $(P_t)_{t\geq 0}$ the semigroup of $X$, which is defined on the space $\BB_b$ of bounded real-valued Borel functions on $\R$ by 
$$\forall f\in \BB_b\ \forall t\geq 0\ \forall x\in\R \quad P_tf(x)=\E_x(f(X_t))=\dint_{\R}f(y)\,p_t(x,dy).$$
\index{$X,(X_t)_{t\geq 0}$: underlying Markov process}
\index{$(P_t)_{t\geq 0}$: semigroup of $X$}
\index{$\BB_b$: set of bounded real-valued Borel functions on $\R$}
We write $\CC_b$ (resp. $\CC_0$) for the set of continuous bounded (resp. vanishing at $\pm\infty$) real-valued functions on $\R$. We make the usual assumption that the filtration $(\FF_t)_{t\geq 0}$ is complete and right-continuous and that the process $X$ takes values in the set $\DD(\R_+,\R)$ of real-valued càdlàg functions on $\R_+$ (see for instance Theorem 19.15 in \cite{kallenberg}). 
\index{$\CC_b$: set of continuous bounded functions on $\R$}
\index{$\CC_0$: set of continuous functions on $\R$ vanishing at $\pm\infty$}
\index{$\DD(I,E)$: set of $E$-valued càdlàg functions on an interval $I$}

Let $\g:(0,+\infty)\cv\R$ be an upper semi-continuous curve, let also $\g_{0}=\limsup_{t\to0+}\g_t\in[-\infty,+\infty)$, and define
$$\tau_\g=\inf\{s\geq 0,X_s<\g_s\}.$$
\index{$\g$: continuous curve (boundary for the selection)}
\index{$\tau_\g$: crossing time of the curve $\g$ for $X$}
With our assumptions, $\tau_\g$ is a stopping time such that, on $\{\tau_{\g}<+\infty\}$, we have $X_{\tau_\g} \leq \g_{\tau_\g}$.

Our main results relate the hydrodynamic limit of the $N$-BMP to the existence of a curve $\g$ satisfying the following assumption:
\begin{enumerate}[label=(\roman*)]
    \item for all $t>0$, $\Q_{\mu_0}(\tau_\g>t)=\e^{-t}.$
\end{enumerate} 
In other words, (i) means that, under $\Q_{\mu_0}$, $\tau_\g$ follows the exponential distribution $\EE(1)$. We write for $r\in\R$
$$U(r,t)=\Q_{\mu_0}(X_t\leq r,\tau_\g>t),$$
so that under assumption (i), $\e^tU(\cdot,t)$ is the c.d.f. of $X$ conditioned on staying above the curve $\g$ up until $t$, that is 
$$\e^tU(r,t)=\Q_{\mu_0}(X_t\leq r|\tau_\g>t).$$
\index{$U(\cdot,t)$: distribution of $X$ on $\{\tau_\g>t\}$}
Moreover, we rely on the following two additional assumptions on $X$:
\begin{enumerate}[label=(\roman*)]
\setcounter{enumi}{1}
    \item $X$ is a Feller process\footnote{By Feller process, we mean that its semigroup is Feller: it maps $\CC_0$ into $\CC_0$ and is strongly continuous (see \textsection III.2 in \cite{revuz-yor}). We will also use that these semigroups map $\CC_b$ into $\CC_b$ (which we can deduce from the former assumption by writing any bounded continuous function as a limit of bounded functions that vanish at infinity). },
    \item $X$ is stochastically non-decreasing, \ie\ for every non-decreasing bounded function $f$, $x\leq y$ and $t\geq 0$, we have $P_tf(x)\leq P_tf(y)$.
\end{enumerate} 

Given an integer $N\geq 2$, we denote by $X^1_t,\cdots,X^N_t$ the positions
 of the particles in the $N$-BMP at time $t\geq 0$, and by $F^N(\cdot,t)$ the empirical c.d.f. of the $N$-BMP:
$$F^N(r,t)=\dfrac{1}{N}\dsum_{i=1}^N\1_{\{X^i_t\leq r\}}.$$
\index{$F^N$: empirical c.d.f. of the $N$-BMP}
We will always consider the case where the initial positions of the particles are i.i.d. with common distribution $\mu_0$, and we denote by $\P$ the corresponding probability measure governing the $N$-BMP. 

\begin{theo}
\label{th:hydro-nbmp}
Consider a time horizon $T>0$, and assume (i),(ii) and (iii) are fulfilled. Then for every $t$ in $[0,T]$, $\n{F^N(\cdot,t)-\e^tU(\cdot,t)}\infty$ goes a.s. to $0$ as $N$ goes to infinity, where $\n{\cdot}\infty$ is the supremum norm on $\BB_b$. Furthermore, for every $\b\in\left(0,\frac12\right)$, there exist positive constants $c_1,c_2>0$ depending only on $\b$ and $T$ such that for all $N$, we have 
\begin{equation}
\label{eq:hydro-nbmp-quant}
    \P\left(\n{F^N(\cdot,t)-\e^tU(\cdot,t)}\infty>N^{-\b}\right)\leq c_2\e^{-c_1N^{\frac{1-2\b}3}}.
\end{equation} 
\end{theo}
\index{$\n{\cdot}\infty$: supremum norm on $\BB_b$}

A natural question is then whether, as $N$ goes to infinity, the trajectory of the lowest particle position in the $N$-BMP approaches $\g$. We investigate this question in \thref{th:min}. To state our result, we need an additional assumption on $X$, known in the literature as \textit{Ray's estimate}: denote by $\tau_y$ the (downwards) crossing time of $y$ by $X$, i.e.
$$ \tau_y=\inf\{s\geq 0, X_s < y\},$$
and assume that:
\begin{enumerate}[label=(\roman*)]
    \setcounter{enumi}{3}
    \item for all real numbers $x>y$ and $\k>0$, we have 
    \begin{equation}
    \label{eq:ray}
        \Q_x(\tau_y\leq t)=o(t^\k),
    \end{equation} 
    \item the curve $\g$ is continuous on $(0,+\infty)$.
\end{enumerate}
\index{$\z^\pm_y$: upwards/downwards crossing time of $y$ for $X$}
Let $m^N_t=\min\limits_{1\leq i\leq N} X^i_t$ be the lowest particle position at time $t$.
\index{$m^N_t$: lowest particle position of the $N$-BMP at time $t$}

\begin{theo}
\label{th:min}
Consider $0<t_0<T$ and assume (i) to (v) are fulfilled. Then, almost surely, the trajectory $t\longmapsto m^N_t$ converges to the curve $t\longmapsto\g_t$, uniformly on $[t_0,T]$ as $N$ goes to infinity. Furthermore, for any $\ee>0$ and $\varkappa>0$ there exists a constant $C>0$ (depending on $\ee,\varkappa,T$ and also on $X$) such that for all $N$, we have 
\begin{equation}\label{eq:min}\P\left(\supl_{t\in[t_0,T]}|m^N_t-\g_t|>\ee\right)\leq CN^{-\varkappa}.\end{equation}
In addition, these results still hold with $t_0=0$ provided that $\g$ is continuous at $0$. 
\end{theo}

\subsection{Discussion}
\label{subsec:bib}

\subsubsection{Assumptions}

Our initial understanding of Assumption (i) was as a probabilistic reformulation of an existence result for a free boundary problem, as discussed in the next paragraph. After (almost) completing the present manuscript, we realized that Assumption (i) had already been studied in the context of the so-called inverse first-passage time problem, for which a substantial body of literature exists; our main reference is Klump and Savov \cite{klump23}, whose introduction contains an overview of the subject and relevant works. In \cite{klump23}, it is proved that Assumption (i) is satisfied for every Feller process such that $\Q_{\mu_0}(X_t \in \cdot)$ is diffuse for all $t$ and such that $\Q_{\mu_0}(\tau_{\g} = \tau'_{\g})=1$, where $\tau'_\g=\inf\{s\geq 0,X_s \leq \g_s\}$; it is also proved that these assumptions are satisfied by a wide class of Lévy processes and diffusion processes.

One can also understand Assumption (i) as equivalent (at least at a heuristic level, we refrain from giving a precise mathematical statement here) to the existence of a solution to the following problem:
\begin{equation}
\label{eq:fbp}
    \left\{\begin{array}{ll}
\d_tu(x,t)=\L u(x,t)+u(x,t), &x>\g_t,\\
u(x,t)=0,&x\leq\g_t,
\end{array}\right.
\end{equation}
where $u(\cdot,t)$ is a probability density function, $\L$ is the infinitesimal generator of $X$, and $\g$ is part of the problem as a free boundary. In the case where $X$ is a standard Brownian motion, so that $\L=\frac12\d^2_x$, this problem has been solved by Berestycki et al. in \cite{fbp-fkpp} and by Lee\footnote{To be precise, Berestycki et al. prove the existence of global solution for \eqref{eq:fbp} with $\L=\frac12\d^2_x$, while Lee proves in \cite{lee} the existence of local solutions to a more general version of the free boundary problem corresponding to the version of the $N$-BBM with non-local branching that appears in \cite{nonlocal-nbbm}.} in \cite{lee}.

Another situation in which we know that Assumption (i) is verified is the case of Markov processes starting from a quasi-stationary distribution on $[0,+\infty)$, in which case $\g$ is constant: it is well-known that in this case, the first crossing time of the origin is exponentially distributed\footnote{See for instance \cite{qsd} Theorem 2.2. Notice that we might need to rescale the process in time to obtain the exponential distribution with parameter $1$. }. Various authors have given some conditions under which such distribution exist, starting with the case of the Brownian motion with negative drift \cite{qsd-bm}, followed by more general one-dimensional diffusion (see \cite{qsd} and references therein). Several conditions have also been given for Lévy processes \cite{kyprianou-palmowski,qsd-levy,yamato}.

Assumption (ii) (Feller property) is  quite classical for Markov processes and we do not discuss it any further. 

Assumption (iii) (stochastic monotonicity) is automatically verified for Feller processes with continuous paths as well as Lévy processes (see \cite{BerFre} for a more detailed discussion about this assumption). 

Assumption (iv) (Ray's estimate) is verified for any Feller process with continuous paths (see Section 3.8 in \cite{ito-mckean}).

Assumption (v) is, to our knowledge, established for standard Brownian motion (see \cite{fbp-fkpp, chen11}) and, more generally, diffusions with smooth coefficients \cite{chen22}).

\subsubsection{Related works}

The hydrodynamic limit shown in \thref{th:hydro-nbmp} has already been established by De Masi et al. for the $N$-BBM in \cite{nbbm,nonlocal-nbbm} (allowing a more general non-local branching mechanism in \cite{nonlocal-nbbm}), with a control upon the probabilities of deviations for the supremum norm given by a negative power of $N$. The multi-dimensional case (the so-called Brownian bees) has also been studied by Berestycki et al. in \cite{bees}, where this time the killed particle is the furthest away from the origin. Taking euclidean norms of particles' positions, this model leads to an $N$-BMP driven by the opposite of a Bessel process. Berestycki et al. not only obtain the hydrodynamic limit, but also a convergence result analogous to our \thref{th:min}, though only at a fixed time $t$ instead of the uniform convergence over $[t_0,T]$. 

As in \cite{nbbm,nonlocal-nbbm,bees}, our proof relies on a comparison between the $N$-BMP and auxiliary stochastic processes. However, our version of these processes is based on a different definition where the curve $\g$ plays a central role. This approach allows us to get improved bounds (stretched exponential vs. negative power) for the hydrodynamic limit, thus allowing us to also prove the uniform convergence of the lowest position, extending the strategy used in \cite{bees}\footnote{For this last problem, the negative power bound we obtain comes from the use of Ray's estimate, which holds for very general Markov processes; in the Brownian case, our strategy would lead to a stretched exponential bound by using the stronger control on \eqref{eq:ray} available in this specific case.}. 
It is worth mentioning that the result obtained in \thref{th:hydro-nbmp} is universal: the constants that control the deviations of the distribution of the $N$-BMP do not depend on the underlying Markov process $X$ or the boundary $\g$ (as soon as the assumptions are verified). However, in \thref{th:min}, the constants depend on $X$ both through Ray's estimate and the moduli of continuity of $\g$.  

In addition to these works, the hydrodynamic limit has also been investigated for a variety of branching-selection systems with a different selection and/or branching mechanism, e.g. \cite{durrett-remenik,groisman20, groisman21, atar23}, often in connection with free boundary problems.

\subsection{Outline of the proof}

To establish the hydrodynamic limit in \cite{nbbm,nonlocal-nbbm,bees}, a key idea is to compare the particle process with stochastic barriers in which selection is either delayed or anticipated so that it occurs at fixed times rather than at the branching times. Here we use a different approach that takes advantage of the assumed existence of a curve $\g$ satisfying condition (i), and compare the $N$-BMP to another process that we call the $\g$-BMP. To be specific, starting from $n$ i.i.d. initial positions with distribution $\mu_0$, the $\g$-BMP is obtained by running a (dyadic) branching Markov process (BMP) based on the dynamics of  $X$, while removing (or "killing") particles as soon as they cross the moving boundary $\g$. When $\g$ satisfies condition (i), the many-to-one formula (see \hyperref[lem:many-to-one]{Lemma \ref{lem:many-to-one}} below) implies that the expected population size of the $\g$-BMP at any given time $t$ is equal to $n$, since the exponential growth due to branching is, on average, exactly offset by the exponential decay due to killing below $\g$ while starting with distribution $\mu_0$. This suggests that the $\g$-BMP starting with $n=N$ particles should in some sense be close to the $N$-BMP, where selection keeps the population size exactly equal to $N$. 

The $\g$-BMP turns out to be more tractable than the $N$-BMP, since the way particles are removed in the $\g$-BMP does not involve the complicated interaction implied by the $N$-BMP selection mechanism\footnote{Such an observation is not new: with $\g$ a suitably chosen linear boundary, it plays a substantial role in the study the long-time behaviour of the $N$-BBM \cite{maillard-nbbm} or the $N-$BRW (branching random walk) \cite{berard-gouere}.}. When $n$ is large, we are able to prove in \propref{prop:hydro-gbmp} that, for the $\g$-BMP at time $t$, with high probability, the population size is indeed of order $n$, and the empirical distribution of particle positions is close to $\Q_{\mu_0}(X_t \in \cdot | \tau_{\g} > t)$ (we actually prove explicit, non-asymptotic bounds).

Our idea is then to bound the $N$-BMP $\X$ between two copies $\X^+, \X^-$ of the $\g$-BMP starting respectively with $n=\lceil N(1+\dd)\rceil$ and $n=\lfloor N(1-\dd)\rfloor$ particles. 
More precisely, in \propref{prop:coupling-rapid-proof}, we build  a coupling between $\X$ and $\X^+, \X^-$ in such a way that, provided that $\#\X^+_s\geq N$ for all $s\leq t$, one has $\X_t\lsto\X^+_t$, and provided that $\#\X^-_s\leq N$ for all $s\leq t$, one has $\X^-_t\lsto\X_t$, where, given two real-valued tuples $x,y$ (possibly with different sizes), $x \lsto y$ means that, $\forall r\in\R\quad \#\{i, x^i\geq r\}\leq\#\{i, y^i\geq r\}$. The stochastic monotonicity assumption (iii) is key to the construction of the coupling.

Then, thanks to our control on the population size of the $\g$-BMP, we are able to prove in \propref{prop:smallevents} that, given a finite time horizon $T$, for large $N$, with high probability, the population size of $\X^+$ (resp. $\X^-$) remains above (resp. below) $N$ over the time-interval $[0,T]$, for suitably small values of $\delta$.  In turn, this allows us to compare the empirical distribution of particles within the $N$-BBM at time $t$ to $\Q_{\mu_0}(X_t \in \cdot | \tau_{\g} > t)$, leading to \thref{th:hydro-nbmp}.

The proof of Theorem \thref{th:min} then combines the result of \thref{th:hydro-nbmp} together with Ray's estimate to control the difference between $\g_t$ and $m^N_t$. Our proof is similar in spirit to that of Berestycki et al. in \cite{bees} (namely Section 4.2 about the proof of their Proposition 1.6). Compared to \cite{bees}, we make the most of our improved bounds to obtain uniform convergence, even though we are using a weaker control on the trajectories of $X$ (given by Ray's estimate).

\subsection{Organization of the paper}

The rest of the paper is organized as follows. We start by studying the $\g$-BMP and the stochastic barriers in \secref{sec:gbmp}, where we prove a hydrodynamic limit result for the $\g$-BMP (\propref{prop:hydro-gbmp}) and then deduce \propref{prop:hydro-stobar}. In this section we also study the size of the $\g$-BMP in order to prove \propref{prop:smallevents}. In \secref{sec:proof-nbmp}, we prove \propref{prop:coupling-rapid-proof} 
and deduce \thref{th:hydro-nbmp}. Finally, in \secref{sec:min} we prove \thref{th:min} thanks to the additional regularity assumption (iv). 


\section{Study of the $\g$-BMP}
\label{sec:gbmp}

\subsection{Hydrodynamic limit of the $\g$-BMP}

Recall the definition of the (dyadic) branching Markov process (BMP), a special case of the class of processes studied in \cite{inw}: we have a finite number of particles on the real line following independent trajectories of the Markov process $X$; each particle branches at rate $1$ independently of the others, giving birth to another particle at the same location, which immediately starts to follow yet another independent trajectory of the Markov process $X$, and so on. The corresponding $\g$-BMP is obtained  by removing particles as soon as they cross the moving boundary $\g$.

To be more precise, let us denote by $\mathfrak{X}$ a BMP starting at time $0$ with a single particle whose location follows the distribution $\mu_0$, and let $\NN^\mathfrak{X}_t$ denote the population size of $\mathfrak{X}$ at time $t$. By identifying the trajectory of a particle before its birth to that of its parent, we can consider, for each $j=1,\ldots,\NN^\mathfrak{X}_t$, the trajectory $(X^j_s)_{s \geq 0}$ of the $j-$th particle\footnote{The particle present at time $0$ is given number $1$; successively born particles are then numbered by birth order. The quantities we consider are in fact independent of which numbering scheme we adopt.} at time $t$, and the corresponding crossing time  $\tau^j_\g$ of the boundary.  The number of particles in the corresponding $\g$-BMP whose locations are above $r$ at time $t$, is then given by 
$$\chi_t(\mathfrak{X},r)=\dsum_{j=1}^{\NN^\mathfrak{X}_t}\1_{\{X^j_t>r,\tau^j_\g>t\}}.$$
\index{$\mathfrak X,\mathfrak X^i$: BMP starting from one particle}
\index{$\chi(\mathfrak X,r)$: number of particles above $r$ for the $\g$-BMP made from $\mathfrak X$}

We now recall the many-to-one lemma, which we state here in the version that will be useful throughout the paper (see Lemma 1 and Section 4.1 in \cite{harris-roberts}).

\begin{lemma}[Many-to-one lemma]
\label{lem:many-to-one}
    Suppose that $\mathfrak{X}=(X^i_t)_{t\geq 0,1\leq i\leq \NN^\mathfrak{X}_t}$ is a BMP starting from one particle with distribution $\mu_0$, and let $f:\DD([0,t],\R)\cv\R$ be a bounded measurable function on the set of càdlàg functions from $[0,t]$ to $\R$, for some positive time $t$. Then we have 
    \begin{equation}
    \label{eq:many-to-one}
        \E\left(\dsum_{i=1}^{\NN^\mathfrak{X}_t}f(X^i_{[0,t]})\right)=\e^t\E_{\mu_0}(f(X_{[0,t]})).
    \end{equation}
\end{lemma}
Defining $G(r,t)$ for $r$ in $\R$ and $t\geq 0$ by
\begin{equation}
\label{eq:comp-cdf-Xcond}
    G(r,t)=1-\e^tU(r,t)=\e^t\Q_{\mu_0}(X_t>r,\tau_\g>t)=\Q_{\mu_0}(X_t>r|\tau_\g>t),
\end{equation}
\index{$G$: complement of the c.d.f. of $X$ conditioned on not crossing $\g$}
we deduce from the many-to-one lemma (\hyperref[lem:many-to-one]{Lemma \ref{lem:many-to-one}}) that, for $r$ in $\R$ and $t\geq 0$, 
\begin{equation}\label{eq:en-moyenne-mto}\E(\chi_t(\mathfrak{X},r))=G(r,t).\end{equation}

Now let $\mathfrak{X}^1,\mathfrak{X}^2,\ldots$ denote an i.i.d. sequence of copies of $\mathfrak{X}$, so that, for all $n \geq 1$, $\mathfrak{X}^1,\ldots, \mathfrak{X}^n$ joined together form a BMP starting with $n$ particles. Accordingly, define the population sizes $\NN^i_t$, the positions $(X^{ij}_t)_{1\leq j\leq \NN^i_t}$ and the crossing times $(\tau^{ij}_\g)_{1\leq j\leq N^i_t}$ corresponding to $\mathfrak{X}^i$, for $i \geq 1$, and set 
$$G^n(r,t)=\dfrac{1}{n}\dsum_{i=1}^n\chi_t(\mathfrak{X}^i,r),$$ so that $G^n$ is the complement of the empirical c.d.f. of a $\g$-BMP starting from $n$ particles. 
\index{$\NN^{\mathfrak X},\NN^i$: size of the BMP}
\index{$G^n$: complement of the empirical c.d.f. of $n$ BMPs}

Given \eqref{eq:en-moyenne-mto} and the fact that the random variables $\chi_t(\mathfrak{X}^1,r),  \chi_t(\mathfrak{X}^2,r), \ldots$ are i.i.d., the law of large numbers shows that $G(r,t)$ is the limit of $G^n(r,t)$ as $n$ goes to infinity, so that the hydrodynamic limit of the $\g$-BMP can be identified as the underlying Markov process $X$ conditioned on not crossing the boundary. 
The following result provides a much more precise control on the difference between $G^n$ and $G$, with quantitative non-asymptotic bounds on the sup-norm $\n{G^n(\cdot,t)-G(\cdot,t)}\infty$.

\begin{prop}
\label{prop:hydro-gbmp}
Let $C(t)=2\e^{2t}-\e^t$ and $c_3=-\ln(1-\e^{-T})>0$. Then for any $\eta>0$, $t\in[0,T]$, $n\geq \mathscr N_0=\mathscr N_0(\eta,t) =\max\left(\dfrac{8C(t)}{\eta^2},2\right)$ and $\a\in\left(0,\frac{1}{2}\right)$, we have the following bound:
\begin{equation}
\label{eq:hydro-gbmp}
    \P \left(\n{G^n(\cdot,t)-G(\cdot,t)}\infty>\eta\right)\leq 16n^2\e^{-\frac{n^{1-2\a}\eta^2}{32}}+32\e^Tn^2\e^{-c_3n^\a}.
\end{equation}
\end{prop}
Invoking the Borel-Cantelli lemma, an immediate corollary of \propref{prop:hydro-gbmp} is that, for all $t \geq 0$, almost surely, $\lim_{n \to +\infty} \n{G^n(\cdot,t)-G(\cdot,t)}\infty = 0$. As another corollary, we have the following bound on the total population size of the $\g$-BMP.

\begin{coro}
\label{coro:size-gbmp}
Let $\mathscr X=(X^{ij}_t,1\leq i\leq n,1\leq j\leq \NN^i_t,0\leq t\leq T)$ be a BMP starting from $n$ particles independently distributed according to $\mu_0$, and let $N_t$ be the population size at time $t$ of the $\g$-BMP made from $\mathscr X$. Then\footnote{The definition of $\mathscr N_0$ is given in \propref{prop:hydro-gbmp}.} for any $\eta>0$, $t\in[0,T]$, $n\geq \mathscr N_0(\eta,t)$ and $\a\in\left(0,\frac{1}{2}\right)$, we have 
\begin{equation}
\label{eq:size-gbmp}
\P\left(\left|N_t-n\right|>n\eta\right)\leq 16n^2\e^{-\frac{n^{1-2\a}\eta^2}{32}}+32\e^Tn^2\e^{-c_3n^\a}.
\end{equation}
\end{coro}

\begin{rem} Note that:
\begin{itemize}
    \item The exponential decay is driven by universal constants, in the sense that they depend on $T$ but not on the initial distribution $\mu_0$ nor on the choice of the Markov process that drives the BMP. 
    \item The fact that a $\g$-BMP starting from $n$ particles is the joining of $n$ independent copies of $\g$-BMPs starting with one particle is key to the proof of the proposition, and this is why we study this process rather than studying the $N$-BMP directly. 
    \item Our proof is an adaptation of the proof of Theorem 12.4 in \cite{devroye}, which is a quantitative version of the Glivenko-Cantelli theorem, using symmetrization ideas introduced by Dudley \cite{dudley} and Pollard \cite{pollard}. Here, we adapted this proof to the context of branching Markov processes. 
\end{itemize}
\end{rem}

\begin{demo}[Proof of \propref{prop:hydro-gbmp}. ]

\noindent \textit{Step 1. Symmetrization by a second load of BMPs. }

 Fix an integer $n$ and $\eta>0$, let $\Tilde{\mathfrak{X}}^1,\cdots,\Tilde{\mathfrak{X}}^n$ be $n$ other independent (between themselves and from the $\mathfrak{X}^i$) copies of $\mathfrak{X}$, and use the tildes in the notations for all random objects linked with them. Also write $\mathscr{X}=({\mathfrak{X}}^1,\cdots,{\mathfrak{X}}^n)$ and $\Tilde{\mathscr{X}}=(\Tilde{\mathfrak{X}}^1,\cdots,\Tilde{\mathfrak{X}}^n)$. Let us first prove that for $n\geq \dfrac{8C(t)}{\eta^2}$, we have
\begin{equation}
\label{eq:sym1}
    \P\left(\n{G^n(\cdot,t)-G(\cdot,t)}\infty>\eta\right)\leq 2\P \left(\n{G^n(\cdot,t)-\Tilde G^n(\cdot,t)}\infty>\dfrac{\eta}{2}\right).
\end{equation}
To this end, define the real-valued random variable $r^*$ by 
$$r^*=\begin{cases}
    \inf\{r\in\R,|G^n(r,t)-G(r,t)|>\eta\}\text{ if the set is not empty},\\0\text{ otherwise}.
\end{cases}$$
Note that $r^*$ is $\mathscr X$-measurable. First we have 
\begin{align*}
    \P\left(\n{G^n(\cdot,t)-\Tilde G^n(\cdot,t)}\infty>\dfrac{\eta}{2}\right)&\geq \P\left(|G^n(r^*,t)-\Tilde G^n(r^*,t)|>\dfrac{\eta}{2}\right)\\
     &\geq \P\bigg(|G^n(r^*,t)-G(r^*,t)|>\eta,\\
     &\qquad\qquad\qquad|\Tilde G^n(r^*,t)-G(r^*,t)|<\dfrac{\eta}{2}\bigg). 
\end{align*}
Then we condition on $\mathscr X$ and use that both $G^n(r^*,t)$ and $G(r^*,t)$ are $\mathscr X$-measurable to get
\begin{align}\label{eq:GK-interm}
    &\P\left(\n{G^n(\cdot,t)-\Tilde G^n(\cdot,t)}\infty>\dfrac{\eta}{2}\right) \nonumber \\
    &\qquad\qquad \geq\E \left(\1_{\{|G^n(r^*,t)-G(r^*,t)|>\eta\}}\P\left(|\Tilde G^n(r^*,t)-G(r^*,t)|<\dfrac{\eta}{2}\,\Big|\,\mathscr X\right)\right).
\end{align}
Assume that
\begin{equation}\label{eq:variante-GK}\P\left(|\Tilde G^n(r^*,t)-G(r^*,t)|<\dfrac{\eta}{2}\,\Big|\,\mathscr X\right)\geq \dfrac{1}{2} \mbox{ a.s.} \end{equation}
From \eqref{eq:GK-interm} and \eqref{eq:variante-GK}, one has 
$\P\left(\n{G^n(\cdot,t)-\Tilde G^n(\cdot,t)}\infty>\dfrac{\eta}{2}\right)  \geq \dfrac{1}{2}  \P \left(|G^n(r^*,t)-G(r^*,t)|>\eta\ \right)$, and, noting that 
     $|G^n(r,t)-G(r,t)|$ is greater than $\eta$ for some $r$ in $\R$ if and only if it is for $r=r^*$, \eqref{eq:sym1} follows. We now prove \eqref{eq:variante-GK}. Since $r^*$ is $\mathscr X$-measurable, while $\Tilde{\mathfrak{X}}^1,\cdots,\Tilde{\mathfrak{X}}^n$ are i.i.d. copies of $\mathfrak{X}$ independent from $\mathscr X$, we have
$$\E (\Tilde G^n(r^*,t)|\mathscr X)=G(r^*,t)$$
and
$$\V(\Tilde G^n(r^*,t)|\mathscr X)=\dfrac{1}{n}\V\left(\dsum_{j=1}^{\NN_t^{\mathfrak{X}}}\1_{\{X^j_t>r^*,\tau^j_\g>t\}}\,\Big|\,\mathscr X\right),$$
where ${\mathfrak{X}}$ is assumed to be independent from $\mathscr X$. We can then bound the sum by $\NN^{\mathfrak{X}}_t$, which follows the geometric distribution $\G(\e^{-t})$, and thus we have 
$$\V(\Tilde G^n(r^*,t)|\mathscr X)\leq \dfrac{1}{n}\E \left((\NN^{\mathfrak{X}}_t)^2\right)=\dfrac{C(t)}{n},$$
where $C(t)=2\e^{2t}-\e^t$. 
By Chebyshev's inequality, we obtain
$$\P\left(|\Tilde G^n(r^*,t)-G(r^*,t)|<\dfrac{\eta}{2}\,\Big|\,\mathscr X\right)\geq 1-\dfrac{4C(t)}{n\eta^2},$$
and under the condition $n\geq \dfrac{8C(t)}{\eta^2}$, the conditional probability is indeed greater than $\dfrac{1}{2}$. 

\noindent \textit{Step 2. Second symmetrization with random signs. }

Let $\ee_1,\cdots,\ee_n$ be i.i.d.\,\,Rademacher random variables independent from $\mathscr X$ and $\Tilde{\mathscr X}$. Since $\chi_t(\mathfrak X^i,r)$ and $\chi_t(\Tilde{\mathfrak X}^i,r)$ have the same distribution and are independent for every real number $r$, $\chi_t(\mathfrak X^i,r)-\chi_t(\Tilde{\mathfrak X}^i,r)$ and $\ee_i \left(\chi_t(\mathfrak X^i,r)-\chi_t(\Tilde{\mathfrak X}^i,r) \right)$ also have the same distribution, hence so do $\n{G^n(\cdot,t)-\Tilde G^n(\cdot,t)}\infty$ and $\n{\dfrac{1}{n}\dsum_{i=1}^n\ee_i \left(\chi_t(\mathfrak X^i,\cdot)-\chi_t(\Tilde{\mathfrak X}^i,\cdot)\right)}\infty$. Then 
$$\P\left(\n{G^n(\cdot,t)-\Tilde G^n(\cdot,t)}\infty>\dfrac{\eta}{2}\right)=\P\left(\n{\dfrac{1}{n}\dsum_{i=1}^n\ee_i \left(\chi_t(\mathfrak X^i,\cdot)-\chi_t(\Tilde{\mathfrak X}^i,\cdot) \right)}\infty>\dfrac{\eta}{2}\right),$$
and we can now use the triangle inequality and the union bound to obtain 
\begin{align*}
    &\P\left(\n{\dfrac{1}{n}\dsum_{i=1}^n\ee_i \left(\chi_t(\mathfrak X^i,\cdot)-\chi_t(\Tilde{\mathfrak X}^i,\cdot) \right)}\infty>\dfrac{\eta}{2}\right)\\
    &\qquad\qquad\leq\P\left(\n{\dfrac{1}{n}\dsum_{i=1}^n\ee_i\chi_t(\mathfrak X^i,\cdot)}\infty>\dfrac{\eta}{4}\right)+\P\left(\n{\dfrac{1}{n}\dsum_{i=1}^n\ee_i\chi_t(\Tilde{\mathfrak X}^i,\cdot)}\infty>\dfrac{\eta}{4}\right)\\
    &\qquad\qquad=2\P\left(\n{\dfrac{1}{n}\dsum_{i=1}^n\ee_i\chi_t(\mathfrak X^i,\cdot)}\infty>\dfrac{\eta}{4}\right)
\end{align*}
and thus 
\begin{equation}
\label{eq:sym2}
    \P\left(\n{G^n(\cdot,t)-\Tilde G^n(\cdot,t)}\infty>\dfrac{\eta}{2}\right)\leq 2\P\left(\n{\dfrac{1}{n}\dsum_{i=1}^n\ee_i\chi_t(\mathfrak X^i,\cdot)}\infty>\dfrac{\eta}{4}\right).
\end{equation}

\noindent \textit{Step 3. Conditioning on $\mathscr X$ and using Hoeffding's inequality. }

We have to bound the r.h.s. in \eqref{eq:sym2}. The main issue is that we have to estimate a supremum over $\R$. However, conditionally on $\mathscr X$, we actually take the supremum over a finite number of random variables (only the $\ee_i$ remain random), which enables us to use a union bound. 

To see this, let us consider $n$ realizations $x^1,\cdots,x^n$ of the BMPs $\mathfrak X^1,\ldots, \mathfrak X^n$, with respective population sizes $n^i_t$ at time $t$, and $n_t=\dsum_{i=1}^nn^i_t$ the total size of the population, so that particles present at time $t$ are indexed by pairs $(i,j)$ with $1 \leq i \leq n$ and $1 \leq j \leq n^i_t$. Ordering particles according to their position at time $t$ leads to an ordering of the pairs of the form $(i_k,j_k)$, for $1 \leq k \leq n_t$, such that  the sequence $(x^{i_kj_k}_t)_{1\leq k\leq n_t}$ is non-decreasing.
 Then, as $r$ runs over the real line, $\dsum_{i=1}^n\ee_i\chi_t(x^i,r)=\dsum_{k=1}^{n_t}\ee_{i_k}\1_{\{x^{i_kj_k}_t>r,\tau^{i_kj_k}_\g>t\}}$ takes at most $n_t+1$ values, because each indicator function can change its value at most once (and it may never do, if the crossing time of the boundary is less than $t$). Thus taking the supremum over $\R$ of the values of $\dsum_{i=1}^n\ee_i\chi_t(x^i,r)$ is in fact taking the maximum over $n_t+1$ random variables. We obtain that the event on the r.h.s. of \eqref{eq:sym2}, conditionally on $\mathscr X$, is a finite union, and we dominate its conditional probability by the number of events in the union times the maximum probability of one of these events, which can again be written as a supremum over $\R$: 
\begin{equation}
\label{eq:sup-bc-cond}
    \P\left(\sup_{r\in\R}\left|\dfrac{1}{n}\dsum_{i=1}^n\ee_i\chi_t(X^i,r)\right|>\dfrac{\eta}{4}\,\Big|\,\mathscr X\right)\leq(\NN_t+1)\supl_{r\in\R}\P\left(\left|\dfrac{1}{n}\dsum_{i=1}^n\ee_i\chi_t(X^i,r)\right|>\dfrac{\eta}{4}\,\Big|\,\mathscr X\right),
\end{equation}
where $\NN_t=\dsum_{i=1}^n\NN^i_t$. Now we use the Hoeffding inequality for sums of independent bounded random variables (see for instance Theorem 8.1 in \cite{devroye}) to bound the probability in the r.h.s. of \eqref{eq:sup-bc-cond}, uniformly with respect to $r$: conditionally on $\mathscr X$, $\ee_i\chi_t(X^i,r)$ is bounded by $\pm \NN^i_t$ and centered, hence, setting $A_t:= \exp\left(-\textstyle{\frac{n^2\eta^2}{32\sum\limits_{i=1}^n(\NN^i_t)^2}} \right)$, we have that
$$\P\left(\left|\dfrac{1}{n}\dsum_{i=1}^n\ee_i\chi_t(X^i,r)\right|>\dfrac{\eta}{4}\,\Big|\,\mathscr X\right)\leq 2 A_t.$$
Now since the bound is uniform with respect to $r$,  \eqref{eq:sup-bc-cond} leads to 
\begin{equation}
\label{eq:sup-bc-at}
\P\left(\sup_{r\in\R}\left|\dfrac{1}{n}\dsum_{i=1}^n\ee_i\chi_t(X^i,r)\right|>\dfrac{\eta}{4}\,\Big|\,\mathscr X\right)\leq2(\NN_t+1)A_t.
\end{equation}

\noindent \textit{Step 4. Final domination.}

Thanks to steps 1 to 3, we obtain by concatenating \eqref{eq:sym1}, \eqref{eq:sym2} and \eqref{eq:sup-bc-at}
\begin{equation}
\label{eq:At0}
\P\left(\n{G^n(\cdot,t)-G(\cdot,t)}\infty>\eta\right)\leq 8\E\left((\NN_t+1)A_t\right).
\end{equation}
Now we fix $K>0$ and split the r.h.s. depending on whether all $\NN^i_t$ are smaller than $K$ or not. 
\begin{itemize}
    \item On $\{\NN_t^i\leq K,\forall i\}$, then $A_t$ is smaller than $\e^{-\frac{n\eta^2}{32K^2}}$ and $\NN_t$ smaller than $nK$, so that 
    \begin{equation}
    \label{eq:At1}
        \E \left((\NN_t+1)A_t\1_{\{\NN_t^i\leq K,\forall i\}}\right)\leq (nK+1)\e^{-\frac{n\eta^2}{32K^2}}.
    \end{equation}
    \item For the complement, we use the union bound and get
    $$\E \left((\NN_t+1)A_t\1_{\{\NN_t^i\leq K,\forall i\}^c}\right)\leq\E \left(\left(\dsum_{i=1}^n\NN^i_t+1\right)A_t\dsum_{i=1}^n\1_{\{\NN_t^i> K\}}\right),$$
    then expand the sums and use independence: the $\NN^i_t$ are i.i.d. random variables with geometric distribution and $A_t$ is smaller than $1$, so we obtain
    \begin{align}
        &\E \left(\left(\dsum_{i=1}^n\NN^i_t+1\right)A_t\dsum_{i=1}^n\1_{\{\NN_t^i> K\}}\right)\notag\\
        &\qquad\qquad\leq n\P(\Gg>K)+n\E (\Gg\1_{\{\Gg>K\}})+n(n-1)\E (\Gg)\P(\Gg>K),
    \end{align}
    where $\Gg$ has geometric distribution with parameter $\e^{-t}$. Finally we use \lemref{lem:geom} in \secref{sec:geom} to get 
    \begin{equation}
    \label{eq:At2}
    \E \left((\NN_t+1)A_t\1_{\{\NN_t^i\leq K,\forall i\}^c}\right)\leq n\e^{-c(t)K}+n(K+2)\e^t\e^{-c(t)K}+n(n-1)\e^t\e^{-c(t)K},
    \end{equation}
    with $c(t)=-\ln(1-\e^{-t})>0$.  
\end{itemize}
We now take $K=n^\a$ with $\a\in\left(0,\frac12\right)$, and use that $c(t)\geq c_3:=c(T)$ and $\e^t\leq\e^T$ to get the following bound (which is uniform with respect to $t$): then \eqref{eq:At2} gives for all $n$
\begin{equation}
\label{eq:At2bis}
    \E \left((\NN_t+1)A_t\1_{\{\NN_t^i\leq K,\forall i\}^c}\right)\leq (n+n(n^\a+2)\e^T+n(n-1)\e^T)\e^{-c_3n^\a}.
\end{equation}
Summing up step 4, we take \eqref{eq:At0}, \eqref{eq:At1} and \eqref{eq:At2bis} together to obtain
$$\P\left(\n{G^n(\cdot,t)-G(\cdot,t)}\infty>\eta\right)\leq 8(n^{1+\a}+1)\e^{-\frac{n^{1-2\a}\eta^2}{32}}+8(n+n(n^\a+2)\e^T+n(n-1)\e^T)\e^{-c_3n^\a}.$$
Now, as soon as $n\geq 2$ (which is implied by the condition on $n$ obtained in step 1 when $\eta$ is small anyway), we can simplify this bound using $n^{1+\a}+1\leq 2n^2$ and $n+n(n^\a+2)\e^T+n(n-1)\e^T\leq 4n^2\e^T$, and we obtain \eqref{eq:hydro-gbmp}. \end{demo}

\begin{demo}[Proof of \cororef{coro:size-gbmp}. ]Notice that almost surely 
$$\liml_{r\to-\infty}G^n(r,t)=\dfrac{N_t}n$$
while $G(r,t)$ goes to $1$ as $r$ goes to $-\infty$, so that 
$$\P\left(\left|N_t-n\right|>n\eta\right)\leq\P\left(\n{G^n(\cdot,t)-G(\cdot,t)}\infty>\eta\right),$$
and we now deduce \eqref{eq:size-gbmp} from \eqref{eq:hydro-gbmp}. 
\end{demo}

\subsection{Application to $\X^+$ and $\X^-$}\label{sec:appli-barriers}

Fix $\dd\in\left(0,\frac12\right)$ and consider two $\g$-BMPs $\X^{\pm,\dd,N}$ respectively starting from $N^+_\dd:=\lceil N(1+\dd)\rceil$ and $N^-_\dd:=\lfloor N(1-\dd)\rfloor$ i.i.d.\,\,positions with distribution $\mu_0$. We will drop the dependence of the processes in $\dd$ and $N$ in the notations from now, and simply write $\X^+$ and $\X^-$. We will refer to these processes as \textit{upper} and \textit{lower stochastic barrier}: we use the same terminology as in the related works even though the construction is different because these processes play a similar role in the proof. 
\index{$\X^\pm$: stochastic barriers}
\index{$N^+_\dd$: $\lceil N(1+\dd)\rceil$, initial size of $\X^+$}
\index{$N^-_\dd$: $\lfloor N(1-\dd)\rfloor$, initial size of $\X^-$}

In \propref{prop:hydro-stobar} below, we specialize the result of \propref{prop:hydro-gbmp} to these stochastic barriers.
Let us define $G^{\pm,\dd,N}(r,t)=\dfrac{1}{N}\dsum_{x\in \X^{\pm}_t}\1_{\{x>r\}}$ (note that we normalize by $N$ instead of $N^{\pm}$). 
 \begin{prop}
\label{prop:hydro-stobar}
 For all $\dd, \a$ in $\left(0,\frac12\right)$ and $\eta$ in $(0,1)$ there exist an explicit integer $\mathscr N_1=\mathscr N_1(\eta,\dd,T)$ and a constant $c_3=c_3(T)>0$ such that for all $t$ in $[0,T]$ and $N\geq \mathscr N_1$, we have
\begin{equation}
\label{eq:hydro-stobar}
\P\left(\n{G^{\pm,\dd,N}(\cdot,t)-(1\pm\dd)G(\cdot,t)}\infty>\eta\right)\leq 256N^2\e^{-\frac{\eta^2\left(\frac12-\dd\right)^{3-2\a}N^{1-2\a}}{128}}+512N^2\e^{-c_3\left(\frac12-\dd\right)^\a N^\a}.
\end{equation}
\end{prop}
\index{$G^{\pm,\dd,N}$: complement of the empirical c.d.f. of $\X^\pm$}

\begin{demo}[Proof of \propref{prop:hydro-stobar}. ]
By definition, we have 
$$\P\left(\n{G^{\pm,\dd,N}(\cdot,t)-(1\pm\dd)G(\cdot,t)}\infty>\eta\right)=\P_{\mu_0}\left(\n{G^{N_\dd^\pm}(\cdot,t)-\textstyle{\frac{N(1\pm\dd)}{N_\dd^\pm}}G(\cdot,t)}\infty>\textstyle{\frac{N}{N_\dd^\pm}}\eta\right),$$
where $G^{N_\dd^\pm}=\textstyle{\frac{N}{N_\dd^\pm}}G^{\pm,\dd,N}$ is the complement of the empirical c.d.f. of a $\g$-BMP starting with $n=N^\pm_\dd$ particles, as in \propref{prop:hydro-gbmp}. Now we have for $N\geq \textstyle{\frac2\eta}$
$$\left|\textstyle{\frac{N(1\pm\dd)}{N^\pm_\dd}}-1\right|\leq \textstyle{\frac{N}{2N^\pm_\dd}}\eta,$$
and since $G$ is bounded by $1$, we obtain
$$\n{\textstyle{\frac{N(1\pm\dd)}{N_\dd^\pm}}G(\cdot,t)-G(\cdot,t)}\infty\leq \textstyle{\frac{N}{2N_\dd^\pm}}\eta,$$
and hence 
$$\P \left(\n{G^{\pm,\dd,N}(\cdot,t)-(1\pm\dd)G(\cdot,t)}\infty>\eta\right)\leq\P \left(\n{G^{N_\dd^\pm}(\cdot,t)-G(\cdot,t)}\infty>\dfrac{N}{2N_\dd^\pm}\eta\right).$$
We will now use \propref{prop:hydro-gbmp}: let $\eta'=\textstyle{\frac{N}{2N^\pm_\dd}}\eta$, so that the condition on $N^\pm_\dd$ is that $N^\pm_\dd\geq\max\left(\textstyle{\frac{8C(T)}{\eta'^2},\frac2\eta,2}\right)$. We can choose $\eta<1$, so that the second condition implies the third one. Now, since $N^\pm_\dd$ is always between $N(1-\dd)-1$ and $N(1+\dd)+1$, the second condition is achieved as soon as $N\geq \textstyle{\frac{\frac2\eta+1}{1-\dd}}.$
Furthermore, for $N\geq 2$, $\textstyle{\frac{N^\pm_\dd}{N}}$ is between $\frac12-\dd$ and $\textstyle{\frac32}+\dd$, so the first condition is 
$$N^\pm_\dd\geq \textstyle{\frac{32C(T)}{\eta^2}} \left(\textstyle{\frac32}+\dd\right)\qquad\ie\qquad N\geq \frac{ \frac{32C(T)}{\eta^2}  \left(\textstyle{\frac32}+\dd\right)+1}{1-\dd}.$$
We can now write 
$$\mathscr N_1=\mathscr N_1(\dd,\eta,T)=\max\left(\dfrac{\frac2\eta+1}{1-\dd},  \frac{ \frac{32C(T)}{\eta^2}  \left(\textstyle{\frac32}+\dd\right)+1}{1-\dd} \right),$$
so that for $N\geq \mathscr N_1$, \eqref{eq:hydro-gbmp} in \propref{prop:hydro-gbmp} gives for $\a$ in $\left(0,\frac12\right)$
$$\P \left(\n{G^{\pm,\dd,N}(\cdot,t)-(1\pm\dd)G(\cdot,t)}\infty>\eta\right)\leq 16(N^\pm_\dd)^2\e^{-\frac{\eta'^2(N^\pm_\dd)^{1-2\a}}{32}}+32(N^\pm_\dd)^2\e^{-c_3(N^\pm_\dd)^\a},$$
which can be rewritten \eqref{eq:hydro-stobar}
$$\P \left(\n{G^{\pm,\dd,N}(\cdot,t)-(1\pm\dd)G(\cdot,t)}\infty>\eta\right)\leq 256N^2\e^{-\frac{\eta^2\left(\frac12-\dd\right)^{3-2\a}N^{1-2\a}}{128}}+512N^2\e^{-c_3\left(\frac12-\dd\right)^\a N^\a}$$
as desired using again the bounds on $\textstyle{\frac{N^\pm_\dd}{N}}$. 
\end{demo}

As explained in the introduction, we expect that, for large $N$, $\X^+$ (resp. $\X^-$) will have more (resp. less) than $N$ particles with high probability on any compact time interval. We thus define the favorable events 
$$\left\{
\begin{array}{c}
    \Gg^+_t=\Gg^{+,\dd,N}_t=\left\{\#\X^+_s\geq N,\forall s\leq t\right\},\\
    \Gg^-_t=\Gg^{-,\dd,N}_t=\left\{\#\X^-_s\leq N,\forall s\leq t\right\}.
\end{array}
\right.
$$
\index{$\Gg_t^\pm$: favorable on which we can build a coupling $\X\lsto\X^+$ (resp. $\X^-\lsto\X$)}
The following proposition provides the control we need on the probability of these events for large $N$, when $\dd$ goes to zero as a negative power of $N$ (this is the case we consider in subsequent proofs\footnote{We chose to display the bounds as in \eqref{eq:smallevent} by taking $\dd=\frac{N^{-\b}}2$ for the sake of clarity, and because this is the form that we need to prove the other results. In fact, the proof leads to the following more explicit although more complicated bound, for any fixed $\dd\in\left(0,\frac12\right)$. For instance we obtain
$$\P_{\mu_0}((\Gg_T^+)^c)\leq 3N\e^{-c_3N^\a}+3N^{1+\a}\left(\e^{-c_3\left(\frac{N\dd}{2}+1\right)}+36N^2\e^{-\frac{\dd^2N^{1-2\a}}{128(1+\dd)^2}}+72\e^TN^2\e^{-c_3N^\a}\right)$$
for any $N$ greater than a certain $\mathscr N_{2,+}(\dd,T)$ explicitly given in the proof of \propref{prop:smallevents}.}).

\begin{prop}
\label{prop:smallevents}
Take any $\b\in\left(0,\frac12\right)$ and set $\dd=\frac{N^{-\b}}2$. Then there exist positive constants $c_4,c_5$ such that for all $N\geq 2$, we have  
\begin{equation}
\label{eq:smallevent}
    \P ((\Gg_T^\pm)^c)\leq c_5\e^{-c_4N^{\frac{1-2\b}3}}.
\end{equation}
\end{prop}

We now prove \propref{prop:smallevents}. Thanks to \cororef{coro:size-gbmp}, we control the size of the stochastic barriers at any fixed time, but the favorable events $\Gg_T^\pm$ depend on the whole dynamic of the barriers between $0$ and $T$, and not only on the position at one given time. The proof of \propref{prop:smallevents} relies on the use of \cororef{coro:size-gbmp} on a finite (though random) set of times. Indeed, to check if $\Gg^\pm_T$ is realized, we only need to check if the cardinalities are big or small enough at the instants when they may change, meaning when a particle appears or disappears. Even better, for $\Gg^+_T$, we only need to check death events: if $N_t\geq N$ for some time $t$, then this will at least stay true until the next death event, because the number of particles will not decrease until then. Similarly, for $\Gg^-_T$, we only check branching events, because $\Gg^-_T$ can only be invalidated at an instant when the number of particles increases. In the following, we make detailed proof for the bound on the probability of $\Gg_T^+$ and quickly explain what happens for $\Gg_T^-$. 

\begin{demo}[Proof of \propref{prop:smallevents}. ]We assume that our $\g$-BMP is obtained by joining $N^+_{\dd}$ i.i.d. $\g$-BMPs obtained from $N^+_{\dd}$ BMPs $\mathfrak{X}^1,\ldots,\mathfrak{X}^{N^+_{\dd}}$ (which we call families) from which particles are removed as soon as they cross the barrier $\g$. We denote by $N^i_t$ the population size at time $t$ of the $i-$th $\g$-BMP, and let $N_t = \sum_{i=1}^{N^+_{\dd}} N^i_t$, so that $\#\X^+_t=N_t$.
\index{$\mathscr X$: BMP starting from $n$ particles}
Then, for $1 \leq i \leq N^+_{\dd}$ and $\ell \geq 1$, define $\s_{i \ell}$ as the $\ell$-th instant of crossing of the $i$-th family.  

Our goal is then to study $\P(N_{\s_{i \ell}}<N, \ \s_{i \ell} \leq T)$ for all $i,\ell$. Let $A_{i \ell} = \{ \s_{i \ell} \leq T \}$. Here we use the fact that the branching families are independent, so we can separate the $i$-th family from the others: by the union bound, we have 
\begin{align}
    \P(N_{\s_{i\ell}}<N,  A_{i \ell} )&\leq \P\left(\left|N_{\s_{i\ell}}-N^+_\dd\right|>N\dd , A_{i \ell} \right)\notag\\
    &\leq\P\left(\left| \sum_{k\neq i}N^{k}_{\s_{i\ell}}-(N^+_\dd-1) \right|>\textstyle{\frac{N\dd}{2}} , A_{i \ell} \right)+\P\left(|N^{i}_{\s_{i\ell}}-1|> \textstyle{\frac{N\dd}{2}} , A_{i \ell} \right).\label{eq:size-gbmp-deathtime-union}
\end{align}
\index{$N_t,N^i_t$: size at time $t$ of the $\g$-BMP made from $\mathscr X,\mathfrak X^i$}
\index{$\s_{i\ell}$: $l$-th instant of death in the $i$-th family of $\mathscr X$}
For the second term in the r.h.s. of \eqref{eq:size-gbmp-deathtime-union}, we note that, on $A_{i \ell}$, one has that $N^{i}_{\s_{i\ell}} \leq \NN^{i}_{T}$,  and that $\NN^{i}_T$ follows a geometric distribution with parameter $e^{-T}$, so that we have for $N\geq \textstyle{\frac2\dd}$
\begin{align}
    \P\left(|N^{i}_{\s_{i\ell}}-1|> \textstyle{\frac{N\dd}{2}} , A_{i \ell} \right)&=\P\left(N^{i}_{\s_{i\ell}}>\textstyle{\frac{N\dd}{2}}+1 ,A_{i \ell}  \right)\notag\\
    &\leq\P\left(\NN^{i}_{T}> \textstyle{\frac{N\dd}{2}}+1\right)\notag\\
    &\leq\e^{-c_3\left(\textstyle{\frac{N\dd}{2}}+1\right)},\label{eq:size-1gbmp}
\end{align}
where $c_3=-\ln(1-\e^{-T})$. For the first term, take $n=N^+_\dd-1$, $\eta=\textstyle{\frac{N\dd}{2}}\left(N^+_\dd-1\right)^{-1}$ and $\eta'=\textstyle{\frac{\dd}{2(1+\dd)}}$, so that
$$n\eta=\textstyle{\frac{N\dd}{2}}\quad\text{and}\quad \eta\geq\eta'.$$
We can thus rewrite 
$$\P\left(\left|\sum_{k\neq i}N^{k}_{\s_{i\ell}}-(N^+_\dd-1)\right|>\textstyle{\frac{N\dd}{2}}, A_{i \ell} \right)=\P\left(\left|\sum_{k\neq i}N^{k}_{\s_{i\ell}}-n\right|>n\eta, A_{i \ell}\right).$$
Then, since for $k\neq i$ the $\s_{i\ell}$ are independent from the $k$-th family's dynamic, we can disintegrate with respect to the value of $\s_{i\ell}$: 
\begin{equation}
\label{eq:size-stobar+-deathtime}
    \P\left(\left|\sum_{k\neq i}N^{k}_{\s_{i\ell}}-n\right|>n\eta, A_{i \ell} \right)=\dint_0^T\P\left(\left|\dsum_{k\neq i}N^{k}_{t}-n\right|>n\eta\right)h_{\s_{i\ell}}(dt),
\end{equation}
where $h_{\s_{i\ell}}$ is the distribution of $\s_{i\ell}$. Now since $\eta\geq\eta'$, we have for all $t$
\begin{equation}
\label{eq:size-gbmp-1}
\P\left(\left|\dsum_{k\neq i}N^{k}_{t}-n\right|>n\eta\right)\leq\P\left(\left|\dsum_{k\neq i}N^{k}_{t}-n\right|>n\eta'\right),
\end{equation}
and this new integrand is exactly the quantity that appears in \cororef{coro:size-gbmp}, which we can apply provided that $N^+_\dd-1$ is greater than $\mathscr N_0(\eta',t)$, for all $t$ in $[0,T]$. Since $\mathscr N_0(\eta',t)$ increases with $t$, we only need to have $N^+_\dd\geq 3$ and $N^+_\dd\geq \dfrac{8C(T)}{\eta'^2}+1$. Thus, if we take $N\geq \mathscr N_{2,+}=\mathscr N_{2,+}(\dd,T)=\max\left(\textstyle{\frac{32(1+\dd)^2C(T)}{\dd^2}}+1,\textstyle{\frac{3}{1+\dd}}\right)$, we can apply \cororef{coro:size-gbmp} and we have 
$$\P\left(\left|\dsum_{k\neq i}N^{k}_{t}-n\right|>n\eta'\right)\leq 16(N^+_\dd-1)^2\e^{-\frac{(N^+_\dd-1)^{1-2\a}\dd^2}{128(1+\dd)^2}}+32\e^T(N^+_\dd-1)^2\e^{-c_3(N^+_\dd-1)^\a}.$$
Now if $N\geq \mathscr N_{2,+}$, we have $N\leq N^+_\dd-1\leq N(1+\dd)$, so we obtain 
$$\P\left(\left|\dsum_{k\neq i}N^{k}_{t}-n\right|>n\eta'\right)\leq16(1+\dd)^2N^2\e^{-\frac{\dd^2N^{1-2\a}}{128(1+\dd)^2}}+32\e^T(1+\dd)^2N^2\e^{-c_3N^\a}.$$
This bound is uniform with respect to $t\in[0,T]$, so we can also use it to bound the l.h.s. of \eqref{eq:size-stobar+-deathtime}. In view of \eqref{eq:size-gbmp-deathtime-union} and \eqref{eq:size-1gbmp}, we have shown for $N\geq \mathscr N_{2,+}(\dd,T)$, that
$$    \P(N_{\s_{i\ell}}<N, A_{i \ell})\leq \e^{-c_3\left(\frac{N\dd}{2}+1\right)}+16(1+\dd)^2N^2\e^{-\frac{\dd^2N^{1-2\a}}{128(1+\dd)^2}}+32\e^T(1+\dd)^2N^2\e^{-c_3N^\a}.$$
This can be simplified for $N$ large enough as 
\begin{equation}
\label{eq:size-gbmp-deathtime}
\P(N_{\s_{i\ell}}<N, A_{i \ell})\leq c_5\e^{-c_4N^{\frac{1-2\b}3}}
\end{equation}
with $c_4,c_5>0$ depending only on $\b$ and $T$ if we recall $\dd=\textstyle{\frac{N^{-\b}}2}$, absorb the polynomials into the decreasing exponentials and take in them the best compromise $\a=\textstyle{\frac{1-2\b}3}$ for the powers of $N$. 

Now since $N_t$ can decrease only when a particle is removed from the $\g$-BMP, and that the time of such a removal must be one of the $\s_{i \ell}$ for some $1 \leq i \leq N_{\dd}^+$ and $\ell \geq 1$, we have the inclusion
$$(\Gg_T^+)^c \subset \bigcup_{1 \leq i \leq  N_{\dd}^+ \atop \ell \geq 1} \{ N_{\s_{i\ell}}<N, A_{i \ell} \}.$$
Noting that the number of crossings in the $i-$th family up to time $t$ is bounded above by $\NN^{i}_T$, we have, for any integer $K \geq 1$, that 
$\{ \NN^{i}_T \leq K \} \cap \bigcup_{\ell \geq 1} \{ N_{\s_{i\ell}}<N, A_{i \ell} \}  \subset   \bigcup_{\ell = 1}^K  \{ N_{\s_{i\ell}}<N, A_{i \ell} \}$. 
By the union bound and \eqref{eq:size-gbmp-deathtime} we thus have that 
$$\P((\Gg_T^+)^c ) \leq \sum_{i=1}^{N_{\dd}^+} \left(  \P(\NN^{i}_T>K)  + \sum_{\ell=1}^K \P(N_{\s_{i\ell}}<N, A_{i \ell}) \right) \leq  N^+_\dd \left( \e^{-c_3K} +Kc_5\e^{-c_4N^{\frac{1-2\b}3}} \right).$$

Finally we take $K= \left\lceil N^{\frac{1-2\b}3} \right\rceil$, use $N^+_\dd\leq N(1+\dd)+1\leq 2(1+\dd)N\leq 3N$ to obtain for all $N\geq 2$ (up to changing the constants $c_4,c_5>0$)
{
\renewcommand\theequation{\ref{eq:smallevent}}
\addtocounter{equation}{-1}
\begin{equation}
    \P((\Gg_T^+)^c)\leq c_5\e^{-c_4N^{\frac{1-2\b}3}}.
\end{equation}
}

As we mentioned earlier, everything works the same way for $\Gg_T^-$, except that we work with branching times instead of death times (which does not modify the argument because we bounded the number of death times by the number of branching times anyway). We let the reader check that the precise bound we could obtain is 
\begin{align}
    \P((\Gg_T^-)^c)&\leq N\e^{-c_3N^\a}+N^{1+\a}\left(\e^{-c_3\left(\frac{N\dd}{2}+1\right)}\phantom{\e^{-\frac{(N(1+\dd)-1)^{1-2\a}\dd^2}{32(1+\dd)^2}}}\right.\notag\\
    &\quad\left.+16N^2\e^{-\frac{(1-\dd)^{-1-2\a}\dd^2N^{1-2\a}}{128\cdot 2^{1-2\a}}}+32\e^TN^2\e^{-c_3\frac{(1-\dd)^\a N^\a}{2^\a}}\right),
\end{align}
true as soon as $N\geq \mathscr N_{2,-}(\dd,T)=\max\left(\textstyle{\frac{32(1-\dd)^2C(T)}{\dd^2}}+1,\textstyle{\frac{4}{1-\dd}}\right)$, which again, up to modifying the constants, leads for all $N\geq 2$ to 
{
\renewcommand\theequation{\ref{eq:smallevent}}
\addtocounter{equation}{-1}
\begin{equation}
    \P((\Gg_T^-)^c)\leq c_5\e^{-c_4N^{\frac{1-2\b}3}}.
\end{equation}
}
\end{demo}

\section{Hydrodynamic limit of the $N$-BMP}
\label{sec:proof-nbmp}

The key ingredient to prove \thref{th:hydro-nbmp} is the definition of a suitable coupling between the $N$-BMP and the $\g$-BMPs $\X^+$ and $\X^-$ studied in \secref{sec:appli-barriers}.
Recall the definitions of the favorable events: 
$$\left\{
\begin{array}{c}
    \Gg^+_T=\left\{\#\X^+_s\geq N,\forall s\leq T\right\}\\
    \Gg^-_T=\left\{\#\X^-_s\leq N,\forall s\leq T\right\}
\end{array}
\right..
$$
and of the ordering $x \lsto y$ between two real-valued tuples $x,y$ (possibly with different sizes): $$\forall r\in\R\quad \#\{i, x^i\geq r\}\leq\#\{i, y^i\geq r\}.$$ 

\begin{prop}
\label{prop:coupling-rapid-proof}
There exists a coupling of an $N$-BMP $\X$\index{$\X$: $N$-BMP} and its upper stochastic barrier $\X^+$ (resp. lower stochastic barrier $\X^-$) such that on $\Gg_T^+$ (resp. $\Gg_T^-$), we have $\X_t\lsto\X^+_t$ (resp. $\X_t^-\lsto\X_t)$ for all $t$ in $[0,T]$. 
\end{prop}


\begin{demo}[Proof of \propref{prop:coupling-rapid-proof}. ]A key ingredient in the proof is Theorem 1 in \cite{BerFre}, which states that, provided that the Markov process $X$ satisfies properties (ii) (Feller semigroup) and (iii) (stochastic monotonicity), we can build a (Feller) Markov process $(X^1, X^2)$ whose coordinates are versions of the process $X$, such that, starting from  $(X^1_0,X^2_0)=(x^1,x^2)$, with $x^1\leq x^2$, one has  $X^1_t\leq X^2_t$ for every $t \geq 0$. We call this process the monotone pair coupling.

We now explain the construction of a coupling between $\X$ and $\X^+$. As long as $\# \X_t^+ \geq N$, each particle in $\X$ will be paired with a particle in $\X^+$ in such a way that the former always lies below the latter, so that $\X_t\lsto\X^+_t$.

\begin{itemize}
    \item Initial positions: Recall that $\X$ (resp.  $\X^+$) starts with $N$ (resp.  $N^+_\dd=\lceil N(1+\dd)\rceil$) particles whose initial positions form an i.i.d.  family with common distribution $\mu_0$. Since $N^+_\dd\geq N$, we can always couple the initial configurations in such a way that every particle in $\X$ is paired with a particle in $\X^+$ having the same initial position. 
    \item Between branching/killing events: As long as no branching/killing event affecting $\X$ or $\X^+$ occurs, paired particles evolve according to the monotone pair coupling, so their ordering is preserved, while the remaining (unpaired) particles evolve independently.
    \item Branching events: Branching events affecting paired particles are synchronized, so that paired particles simultaneously branch at rate $1$. Such an event leads to two newborn particles, which are then paired together: this is consistent with the ordering constraint since the positions of the paired parent particles are ordered at the time of birth. Note that the $\X^+$ particle that was paired with the lowest particle in $\X$ is now unpaired, because the latter has just been killed as a consequence of the branching event, but it is still true that each particle in $\X$ is paired with a particle in $\X^+$. Finally, unpaired particles branch independently. 
    \item Crossing events: Remember that particles in $\X^+$ are killed when they cross the moving boundary $\g$. As a consequence, when a paired particle in $\X^+$ crosses the boundary, the particle in $\X$ it was paired to must be given another partner. This is possible as long as $\X^+$ contains at least $N$ particles at the time of the crossing: by definition, each of these particles lies above the barrier, among which one at least is not yet paired, while the particle in $\X$ that was paired to the particle that has just crossed, has to lie below the barrier; as a consequence, we can indeed find a partner in $\X^+$ such that the ordering is preserved. In the case where $\X^+$ contains strictly less than $N$ particles at the time of the crossing, we cease to consider pairs of particles and let the processes $\X$ and $\X^+$ evolve independently.  
\end{itemize}

We can make a similar construction of a coupling between $\X$ and $\X^-$. This time, we ensure that, as long as $\# \X_t^- \leq N$, each particle in $\X^-$ is paired with a particle in $\X$ in such a way that the former always lie below the latter.

\begin{itemize}
    \item Initial positions: We start the $N$-BMP from $N$ initial positions  independently distributed according to $\mu_0$, $N^-_\dd$ of which are shared between $\X$ and $\X^-$.
    \item Between branching/killing events: Paired particles subsequently evolve using the monotone pair coupling, while unpaired particles evolve independently.
    \item Crossing events: When a paired particle in $\X^-$ crosses the boundary, the particle in $\X$ it was paired to stops being paired and starts evolving independently, and it is still true that each particle in $\X^-$ is paired to a particle in $\X$. 
       \item Branching events: Branching events affecting paired particles are synchronized, as in the previous coupling between $\X$ and $\X^+$, and the two newborn particles are paired accordingly. In the case where a particle in $\X^-$ was paired with the lowest particle in $\X$, we have to find it a new partner, since this particle is removed as a consequence of the branching event. This is possible as long as $\X^-$ contains at most $N$ particles at the time of the branching: in this case, at most $N-1$ particles in $\X^-$ are paired with particles in $\X$ just before the branching, so at least one particle in $\X$ is not yet paired. This particle lies above the position of the (now defunct) lowest particle, while the particle in $\X^-$ that was paired to the lowest particle lies below, so we can pair these particles together while preserving the order. In the case where $\X^-$ contains strictly more than $N$ particles, we  cease to consider pairs of particles and let the processes $\X$ and $\X^-$ evolve independently.    
 \end{itemize}

Note that, since $\tau^{\g}$ has the exponential $\EE(1)$ distribution with respect to $\Q_{\mu_0}$,  $\P-$a.s.  branching and crossing events always occur at distinct times, so we do not have to consider such situations in the definition of the coupling.

\end{demo}

Now we can finally prove \thref{th:hydro-nbmp} by putting the different pieces together: on $\Gg_T^\pm$ -- which has a very high probability -- the coupling ensures an ordering of the empirical c.d.f., so that $F^N(\cdot,t)$ is squeezed between two functions that converge to small perturbations of $\e^tU(\cdot,t)$, again up to a small probability. Then taking the perturbation small enough, we obtain the limit for $F^N(\cdot,t)$.

\begin{demo}[Proof of \thref{th:hydro-nbmp}.]
Fix $\b>0$, let $\eta=N^{-\b}$ and focus on
$$\P\left(\n{F^N(\cdot,t)-\e^tU(\cdot,t)}\infty>\eta\right)=\P\left(\n{1-F^N(\cdot,t)-G(\cdot,t)}\infty>\eta\right).$$
We first separate this event in such a way that we will be able to use the couplings: we have
\begin{align*}
    \P\left(\n{F^N(\cdot,t)-\e^tU(\cdot,t)}\infty>\eta\right)&\leq\P\left(\exists r\in\R, 1-F^N(r,t)-G(r,t)>\eta\right)\\
    &\quad\quad+\P\left(\exists r\in\R, 1-F^N(r,t)-G(r,t)<-\eta\right).
\end{align*}
Now for every $r\in\R,\dd>0$ and $\ee>0$, the stochastic ordering $\X_t\lsto\X^+_t$ implies that, on $\Gg^{+,\dd,N}_T$,
$$1-F^N(r,t)\leq G^{+,\dd,N}(r,t),$$
so that if $\n{G^{+,\dd,N}(\cdot,t)-(1+\dd)G(\cdot,t)}\infty\leq\ee$, $1-F^N(r,t)\leq G(r,t)(1+\dd)+\ee$,
which we can rewrite as $1-F^N(r,t)-G(r,t)\leq \dd G(r,t)+\ee\leq \dd+\ee$.
Hence if we take $\dd=\ee=\dfrac{\eta}{2}=\frac{N^{-\b}}2$, 
\begin{align*}
    &\P\left(\exists r\in\R, 1-F^N(r,t)-G(r,t)>\eta\right)\\
    &\quad\quad\quad\quad\quad\quad\quad\quad\leq\P\left((\Gg^{+,\dd,N}_T)^c\right)+\P\left(\left\{\exists r\in\R, 1-F^N(r,t)-G(r,t)>\eta\right\}\cap\Gg_T^{+,\dd,N}\right)\\
    &\quad\quad\quad\quad\quad\quad\quad\quad\leq\P\left((\Gg^{+,\dd,N}_T)^c\right)+\P\left(\n{G^{+,\dd,N}(\cdot,t)-(1+\dd)G(\cdot,t)}\infty>\dfrac{\eta}{2}\right),
\end{align*}
and we can now use the results of Propositions \ref{prop:hydro-stobar} and \ref{prop:smallevents} to bound this quantity: we have for all $\a,\b$ in $\left(0,\frac12\right)$ and $N$ large enough
\begin{align*}
    &\P\left(\exists r\in\R, 1-F^N(r,t)-G(r,t)>\eta\right)\\
    &\qquad\qquad\qquad\leq c_5\e^{-c_4N^{\frac{1-2\b}3}}+256N^2\e^{-\frac{\eta^2\left(\frac{1-\eta}2\right)^{3-2\a}N^{1-2\a}}{512}}+512N^2\e^{-c_1\left(\frac{1-\eta}2\right)^\a N^\a},
\end{align*}
and we can again optimize the value of $\a$ to obtain for some constants $c_1,c_2>0$ and all $N$
\begin{equation}
    \P\left(\exists r\in\R,1-F^N(r,t)-G(r,t)>N^{-\b}\right)\leq \frac{c_2}2\e^{-c_1N^{\frac{1-2\b}3}},
\end{equation}

The same reasoning with the other term, using the coupling with $\X^-$, gives 
\begin{align*}
    &\P\left(\exists r\in\R, 1-F^N(r,t)-G(r,t)<-\eta\right)\\
    &\quad\quad\quad\quad\quad\quad\quad\quad\leq\P\left((\Gg^{-,\dd,N}_T)^c\right)+\P\left(\n{G^{-,\dd,N}(\cdot,t)-(1-\dd)G(\cdot,t)}\infty>\dfrac{\eta}{2}\right),
\end{align*}
which is again bounded by exponentially decaying quantities thanks to Propositions \ref{prop:hydro-stobar} and \ref{prop:smallevents}. Putting these equations together we deduce that 
$$    \P\left(\n{F^N(\cdot,t)-\e^tU(\cdot,t)}\infty>N^{-\b}\right)\leq c_2\e^{-c_1N^{\frac{1-2\b}3}},
$$
which is exactly \eqref{eq:hydro-nbmp-quant}.
 Thanks to the Borel-Cantelli lemma, the almost sure convergence of $F^N(\cdot,t)$ to $\e^tU(\cdot,t)$ in the supremum norm is a straightforward consequence of \eqref{eq:hydro-nbmp-quant}. 
\end{demo}

\section{Convergence of the minimum of the $N$-BMP}
\label{sec:min}

\begin{figure}
    \centering
    \includegraphics[width=\textwidth]{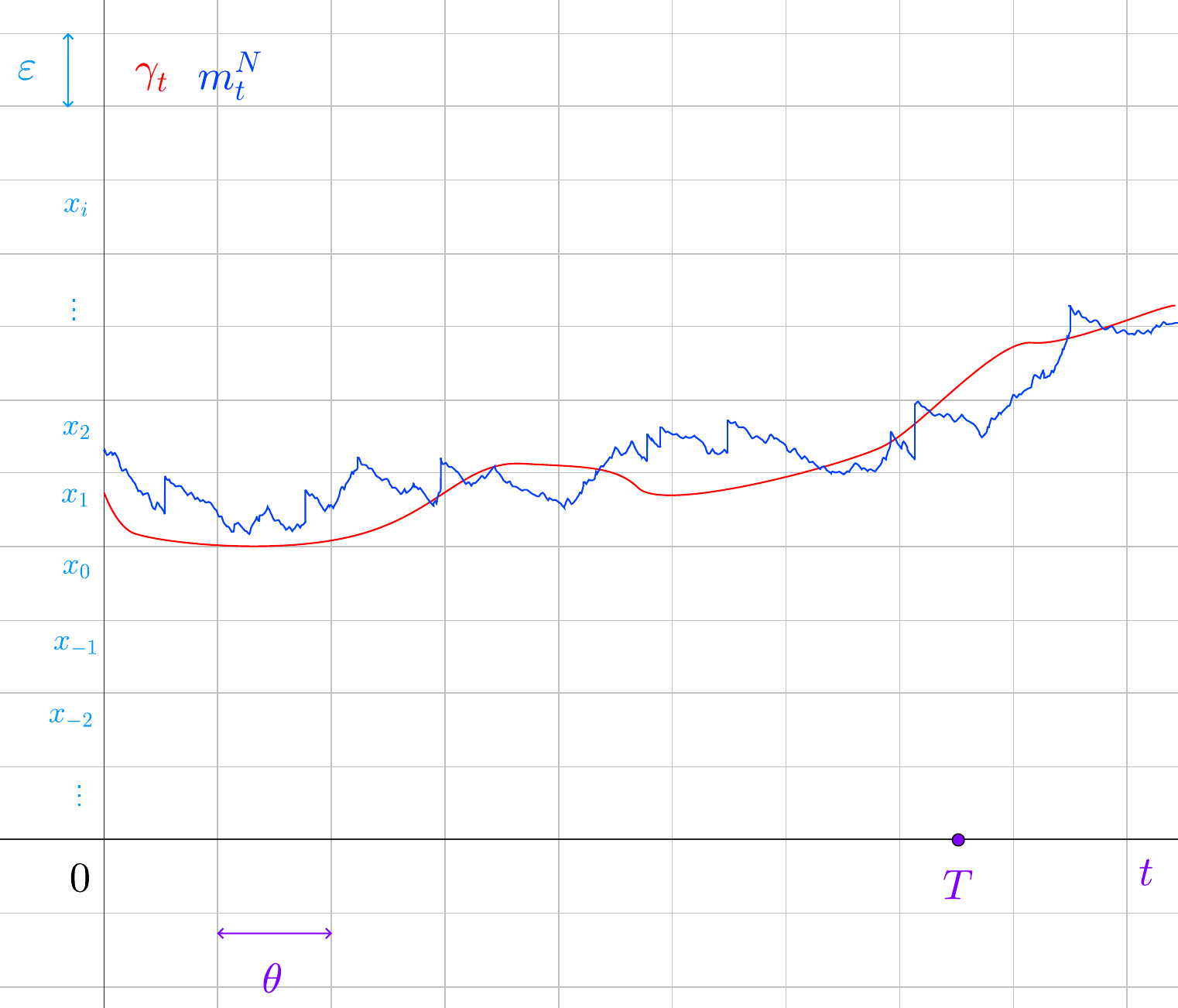}
    \caption{A realization of $m^N_t$ on the space-time grid described in the proof of \thref{th:min}}
    \label{fig:grid-min}
\end{figure}

We now turn to the proof of \thref{th:min} about the convergence of $m^N_t$, which is the lowest particle position of the $N$-BMP at time $t$, towards $\g_t$, when $N$ goes to infinity. As mentioned in the introduction, our proof is similar in spirit to that of Berestycki et al. in \cite{bees} (namely Section 4.2 about the proof of their Proposition 1.6). The idea is to combine the quantitative bound \eqref{eq:hydro-nbmp-quant} of \thref{th:hydro-nbmp} together with Ray's estimate, in order to control the difference between $\g_t$ and $m^N_t$ on a space-time grid (see \figref{fig:grid-min}). Berestycki et al. combine a precise control on Brownian trajectories together with their hydrodynamic limit estimate (involving negative powers of $N$) to prove convergence at any fixed time, while we exploit our improved hydrodynamic limit bound combined with a weaker but more generally valid control on the trajectories (provided by Ray's estimate) to obtain uniform convergence on a time interval. 

\begin{demo}[Proof of \thref{th:min}. ]
We first give the proof assuming that $t_0 > 0$, and then explain the (minor) additional steps needed to deal with the case $t_0=0$.

\noindent \textit{Step 1. Building a suitable space-time grid.}

Fix $\ee>0$ and discretize space with step $\ee$, starting from $\min\limits_{[t_0,T]}\g$, \ie\ set for $i\in\Z$ 
$$x_i=\min\limits_{[t_0,T]}\g+i\ee.$$ 
Now discretize time with step $\t=2N^{-\b}$, for a given $\b\in\left(0,\frac12\right)$, starting from $0$. Let $s_0 \in (0,t_0]$. Since $\g$ is uniformly continuous on $[s_0,T]$, we have, for all large enough $N$ (depending on $\ee,\g, s_0, T, \beta$) that $|\g_t-\g_s|<\ee$ for all $s,t\in[s_0,T]$ such that $|t-s|\leq\t$, and, from now on, we assume that this assumption holds.  Also, the boundedness of $\g$ on $[t_0,T]$ allows us to denote by $K>0$ the smallest integer such that $x_K\geq \max\limits_{[t_0,T]}\g$. 

Thanks to Ray's estimate (Assumption (iv)), for every $\k>0$, there is a constant $C_{\k}$ (depending also on $\ee,t_0,T,\g$) such that\footnote{It is crucial to note that we invoke Ray's estimate for $\Q_{x_{i+1}}(\tau_{x_i}\leq \t)$ simultaneously over a {\bf finite} set of indices $i$.} for all $-6\leq i\leq K+5$, and all $N \geq 1$, 
\begin{equation}
\label{eq:ray-grid}
    \Q_{x_{i+1}}(\tau_{x_i}\leq \t)\leq C_{\k}\t^\k.
\end{equation}
An useful observation is that, due to stochastic monotonicity, this also gives a bound on the probability for one particle to cross the interval $[x_i,x_{i+1}]$ downwards in a short instant:
\begin{equation}
\label{eq:ray-cross}
  \forall x \geq x_{i+1}, \  \Q_{x}(\tau_{x_i}\leq \t)\leq C_{\k}\t^\k.
\end{equation}
To deduce \eqref{eq:ray-cross} from \eqref{eq:ray-grid}, we invoke the monotone pair coupling already used in the proof of \propref{prop:coupling-rapid-proof}, to show that $\Q_{x}(\tau_{x_i}\leq \t) \leq \Q_{x_{i+1}}(\tau_{x_i}\leq \t)$ when $x \geq x_{i+1}$, and \eqref{eq:ray-cross} ensues. Now, using the strong Markov property of $X$, it turns out that \eqref{eq:ray-cross} can be strengthened into
\begin{equation}
\label{eq:ray-cross-strong}
  \forall x, \  \Q_{x}( \exists s \in [0,\theta] \mbox{ such that } X_s \geq x_{i+1} , \  \tau_{x_i}\leq \t)\leq C_{\k}\t^\k.
\end{equation}

\noindent \textit{Step 2. Lower bound at a fixed time. }

Our goal is to establish \eqref{eq:borne-t-fixe}, showing that $\g_t$ is an approximate lower bound for $m^N_t$ at every $t=k\t$ in $(t_0-\t,T]$, where $k$ is an integer. Crucially, the various bounds we derive must hold uniformly with respect to $k$, so we systematically emphasize the dependency of the various constants with respect to the model and grid parameters in the sequel.

We assume that $N$ is large enough (depending on $s_0, t_0, \beta$) so that $t_0-2 \t \geq s_0$, and consider an integer $k$ such that $t:=k\t\in(t_0-\t,T]$. For $i \in \Z$, introduce the numbers
$$r_i=\max\{x_j \mbox{ such that } x_j<\g_t\}-i\ee,$$
and note that $r_0,\ldots, r_5$ are elements of the space grid subject to \eqref{eq:ray-grid}. By definition we have $\g_t>r_0$, and, since $\t$ is an $\ee-$modulus of continuity for $\g$ on $[s_0,T]$ and $s_0 \leq t-\theta \leq t \leq T$, we have that $\g_{t-\t}>r_0-\ee=r_1$.

Now let $\eta=\frac\t2=N^{-\b}$, and, before delving into the details, consider the following informal description of the behaviour of the $N$-BMP over the time-interval $[t-\t,t]$ (see \figref{fig:proof-min}).

\begin{itemize}
    \item Since $\g_{t-\t}>r_1$, the hydrodynamic limit of the $N$-BMP \eqref{eq:hydro-nbmp-quant} implies that with a high probability, fewer than $N\eta$ particles lie below $r_1$ at time $t-\t$.
    \item Comparison with a BMP (without killing) then shows that, with high probability, between time $t-\t$ and $t$, the total number of particles in the $N$-BMP that are, or descend from, particles lying below $r_1$ at time $t-\t$, does not exceed $\sim N\eta\e^\t$. 
    \item With high probability, between $t-\t$ and $t$ we have $\sim N(\e^\t-1)$ newborn particles in the $N$-BMP, so we kill $\sim N(\e^\t-1)$ particles.
    \item Thanks to Ray's estimate, with high probability, particles that did not lie below $r_1$ at time $t-\t$, do not lie below $r_2$ at any time between $t-\t$ and $t$.
    \item As long as there are particles lying below $r_2$, the particles that are killed are among them (because we kill the lowest one).
    \item Since $N(\e^\t-1) \sim 2 N\eta\e^\t$, we deduce from the previous observations that, with high probability, there has to be a time within the interval $[t-\t,t]$ at which no particle is lying below $r_2$.
    \item Hence, thanks to Ray's estimate, with high probability there are no particles lying below $r_3$ at time $t$. 
\end{itemize}

\begin{figure}[h]
    \centering
    \includegraphics[width=\textwidth]{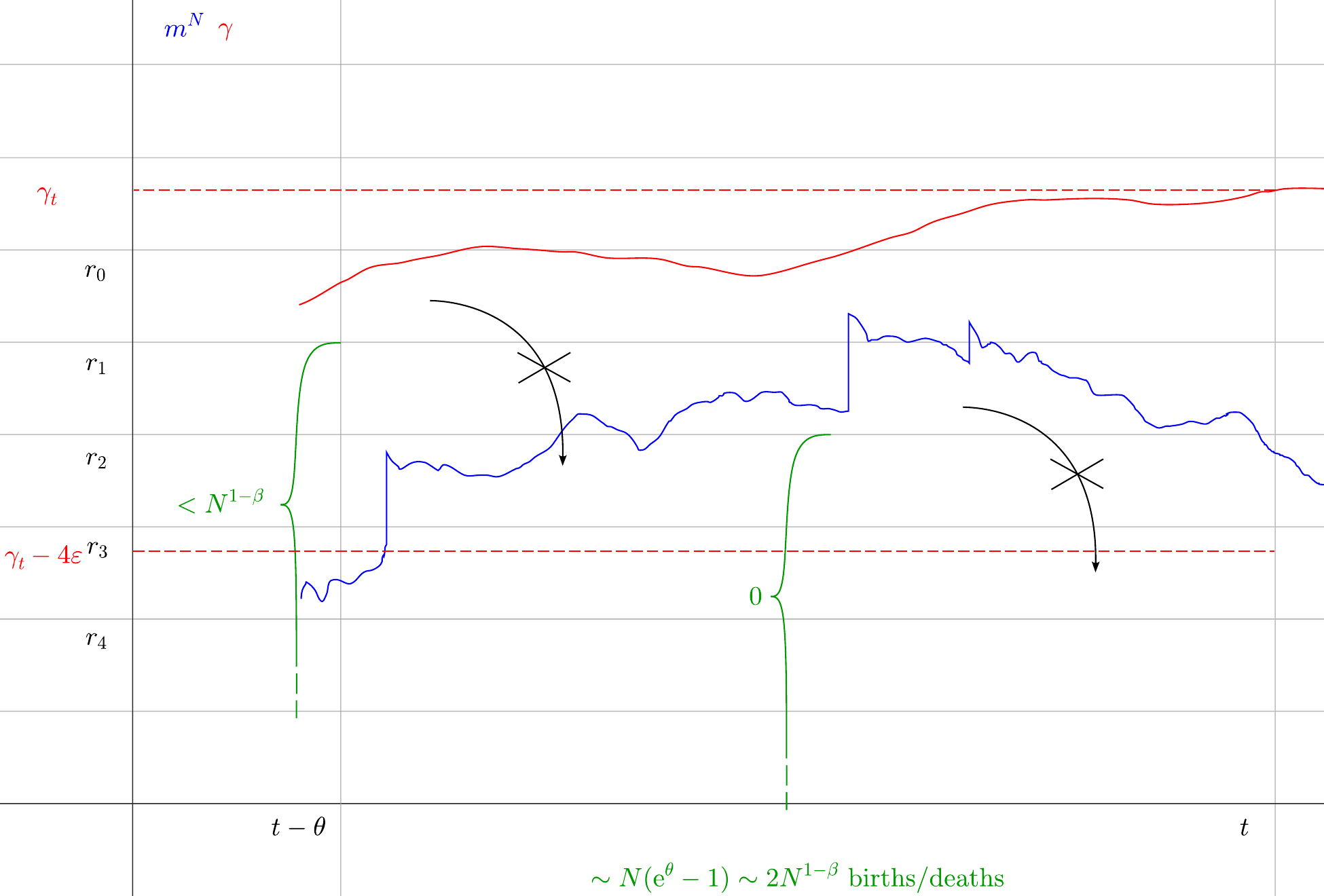}
    \caption{Sketch of the second step of the proof of \thref{th:min}}
    \label{fig:proof-min}
\end{figure}

We now turn the previous informal description into precise arguments. Since $\g_{t-\t}>r_1$, we have that $U(r_1,t-\t)=0$, so that, by \eqref{eq:hydro-nbmp-quant}, 
$$\P\left(|F^N(r_1,t-\t)|>N^{-\b}\right)\leq c_2\e^{-c_1N^{\frac{1-2\b}3}}.$$
Denoting by $K_-$ the number or particles lying below $r_1$ at time $t-\t$, the above bound can be rewritten as
\begin{equation}
\label{eq:pop-below-r1}
    \P(K_->N^{1-\b})\leq c_2\e^{-c_1N^{\frac{1-2\b}{3}}}.
\end{equation}

Now assume that, starting from time $t-\t$, our $N$-BMP is built from a BMP (without killing), which is itself obtained by joining $N$ independent BMPs (or "families"), each starting from a single particle. The $N$-BMP is obtained by  removing particles (together with their descendants) from the BMP, according to the selection mechanism of the $N$-BMP.  

 The population of the BMP at time $t$, hereafter denoted $\NN_t$, is the sum of $N$ independent random variables $G_i$, each with geometric distribution $\G(\e^{-\t})$. For $s$ in $[t-\t,t]$, denote by $X^{ij}_s$ the position of the $j$-th particle from the $i$-th family at time $s$, where $1 \leq j \leq G_i$. For $x>y$, let $\Xi_{x,y}$ denote the event that at least one particle crosses downwards the interval $[y,x]$ on the time interval $[t-\t,t]$, i.e. 
$$\Xi_{x,y}=\{\exists i\in\ent{1}{N}\ \exists j\in\ent{1}{G_i}, (X^{ij}_{t-\t}\geq x\text{ and }\exists s\in[t-\t,t]\ X^{ij}_s< y)\}.$$
\index{$\Xi_{x,y}$: event on which no particle of $\X$ crosses downwards the interval $[y,x]$}
Denote by $\NN'$ the population size at time $t$ of the set of particles that lied below $r_1$ at time $t-\t$ together with their descendants in the BMP, so that $\NN'$ is the sum of the $K^-$ corresponding random variables $G_i$.
Letting $A = \{\forall s\in[t-\t,t], m^N_s < r_2\}$, we now claim that
\begin{equation}\label{eq:inclusion-utile} A \cap \Xi_{r_1,r_2}^c \subset \{ \NN_t-N<\NN' \}.\end{equation}
Indeed, from the condition that the minimum particle position in the $N$-BMP is always $< r_2$ during the time interval $[t-\t,t]$, we deduce that any particle of the $N$-BMP that is killed during this time interval has to lie below $r_2$ at the time of its killing. From the condition that no BMP particle lying above $r_1$ at time $t-\theta$ or one of its descendants, ever goes below $r_2$ during the time interval $[t-\t,t]$, we conclude that every particle that is killed is taken among particles that were below $r_1$ at time $t-\t$ and their descendants in the BMP. Now, the population of particles in the BMP at time $t$ is divided into the $N$ particles that survived the selection mechanism, and the $\NN_t-N$ particles (including their descendants) that have been killed.  We have just seen that this last subset is contained into the subset of $\NN'$ particles (including their descendants in the BMP) that start below $r_1$ at time $t-\t$. Moreover,  the inclusion has to be strict since there is at least one particle below $r_2$ in the $N$-BMP at time $t$, which must be, or descend from, a particle that started below $r_1$ at time $t-\t$, and has survived selection. As a consequence, \eqref{eq:inclusion-utile} is proved.    

Now, using \eqref{eq:inclusion-utile}, we have that  
$$\P(A)=\P(A\cap\Xi_{r_1,r_2})+\P(A\cap\Xi_{r_1,r_2}^c) \leq \P(\Xi_{r_1,r_2})+\P(\NN_t-N<\NN').$$
We note that, conditional upon the value of $K_-$, the random variable $\NN_t-\NN'$ is the sum of $N-K_-$ i.i.d. $\G(\e^{-\t})$ random variables. On the event $\{K_-\leq N^{1-\b}\}$, $\NN_t-\NN'$ is the sum of at least $N-\lfloor N^{1-\b} \rfloor$ such random variables, so that we have 
\begin{equation}
\label{eq:survival_event}
    \P(A)\leq\P(\Xi_{r_1,r_2})+\P\left(\dsum_{i=\lfloor N^{1-\b} \rfloor+1}^NG_i<N\right)+\P(K_->N^{1-\b}).
\end{equation}
We will now deal separately with each of these three terms. 
\begin{itemize}
    \item The rightmost term is dealt with using \eqref{eq:pop-below-r1}.
    \item The leftmost term is bounded thanks to Ray's estimate combined with a large deviations estimate for $\NN_t$, which is a sum of $N$ independent geometric $\G(\e^{-\t})$ random variables. Using the fact that $e^{-\t} = \e^{-2 N^{-\b}} \geq \e^{-2}$, $\NN_t$ is stochastically dominated by a sum of $N$ independent  $\G(\e^{-2})$ random variables, and, since $2 < 1/\e^{-2}$, \lemref{lem:geom} in \secref{sec:geom} ensures that for all $N\geq 1$, $\P(\NN_t>2N)\leq \e^{-c_6N}$, with $c_6=  \Ll^*_{\e^{-2}}(2) > 0$.  Using this bound, we can write 
    $$\P(\Xi_{r_1,r_2})\leq\P(\Xi_{r_1,r_2},\NN_t\leq 2N)+\e^{-c_6N},$$
    and we use a union bound on the first term: we have to deal with at most $2N$ particle trajectories, each of which has a probability smaller than $C_{\k}\t^\k$ to cross the interval $[r_2,r_1]$, thanks to \eqref{eq:ray-cross}, where the value of $\k > 0$ can be taken as large as we want. As a consequence (using also the fact that $\t = 2 N^{-\b}$), we obtain the bound
      \begin{equation}\label{eq:bound-probXi}\P(\Xi_{r_1,r_2})\leq 2NC_{\k}\t^\k+\e^{-c_6N}= 4 C_{\k} N^{1-\b\k}+\e^{-c_6N}.\end{equation}
   
    \item Finally we estimate the central term in \eqref{eq:survival_event}, $\P\left(\dsum_{i=\lfloor N^{1-\b} \rfloor+1}^NG_i<N\right)$, which is again a large deviation probability for a sum of i.i.d. $\G(\e^{-\t})$ random variables. Indeed we have 
    $$\P\left(\sum_{i=\lfloor N^{1-\b} \rfloor+1}^NG_i<N\right)=\P\left(\frac{1}{N-\lfloor N^{1-\b} \rfloor} \sum_{i=\lfloor N^{1-\b} \rfloor+1}^N G_i< \frac{1}{1-N^{-\b}}\right).$$
    For $\t$ small enough, and thus $N$ large enough (depending on $\beta$),
    $\textstyle{\frac{1}{1-\frac\t2}}<\e^\t,$
    so we can indeed use Cramér's theorem to study the large deviations below the mean for this sum of random variables: \lemref{lem:geom} (iv) ensures that
    \begin{equation}
    \label{eq:dev-geom}
        \P\left(\dsum_{i=\lfloor N^{1-\b} \rfloor+1}^NG_i<N\right)\leq\e^{-(N-N^{1-\b})\Ll_{\e^{-\t}}^*\left(\frac{1}{1-N^{-\b}}\right)}.
    \end{equation}
    Since $N^{-\b}=\theta/2$, the expansion \eqref{eq:expansion-leg} yields that $\Ll_{\e^{-\t}}^*\left(\textstyle{\frac1{1-N^{-\b}}}\right)=(1-\ln(2))N^{-\b} +o(N^{-\b})$,
    so that, choosing a constant $c_7$ such that $0<c_{7}<1-\ln(2)$,   \eqref{eq:dev-geom} implies that, for $N$ large enough (depending on $\beta$):
    $$\P\left(\dsum_{i=\lfloor N^{1-\b} \rfloor+1}^NG_i<N\right)\leq \e^{-c_{7}N^{1-\b}}.$$
\end{itemize}
Looking back at \eqref{eq:survival_event}, we have obtained the bound, valid for $N$ large enough (depending on $\b$):
\begin{equation}\label{eq:borne-prob-de-A}\P(A)\leq  4 C_{\k} N^{1-\b\k}+\e^{-c_6N}+\e^{-c_{7}N^{1-\b}}+c_2\e^{-c_1N^{\frac{1-2\b}{3}}}.\end{equation}

Now observe that, on the event $A^c$, there must be a time between $t-\t$ and $t$ at which no particle in the $N$-BMP lies below $r_2$. If, in addition, between $t-\t$ and $t$, no particle is allowed to go below $r_3$ once it has gone above $r_2$, one must then have that $m^N_t\geq r_3$. Thus, using \eqref{eq:ray-cross-strong} and arguing in exactly the same way as we did to obtain the bound \eqref{eq:bound-probXi} on the probability of $\P(\Xi_{r_1,r_2})$, we deduce that, for all large enough $N$ (depending on $\beta$), we have that 
\begin{equation}\label{eq:encore-une-borne}\P(A^c \cap \{m^N_t < r_3 )\}) \leq  4 C_{\k} N^{1-\b\k}+\e^{-c_6N}.\end{equation}
Using the fact that $\g_t-4\ee< r_3$, we have
 $$\P(m_t^N<\g_t-4\ee) \leq \P(m_t^N<r_3) \leq \P(A^c \cap \{m^N_t < r_3 )\})+ \P(A),$$ and combining \eqref{eq:encore-une-borne} and \eqref{eq:borne-prob-de-A}, in which the value of $\kappa$ can be taken as large as we want, we deduce that, for every $\varkappa > 0$, there exists a constant $C^{(1)}$ (depending on $\b,\varkappa,\ee,\g,s_0,t_0,T$, but not on the integer $k$ for which $t=k \t$) such that, for all $N \geq 1$ , 
 \begin{equation}\label{eq:borne-t-fixe}\P(m_t^N<\g_t-4\ee) \leq C^{(1)} N^{-\varkappa}.\end{equation}

\noindent \textit{Step 3. Lower bound on the whole time interval. }

We will use the regularity of $\g$ and $m^N$ to deduce the same kind of bound as \eqref{eq:borne-t-fixe}, but on whole time intervals. Indeed, even though $m^N$ can have some upward jumps (when the particle at the lowest position is killed while being far away from all the other particles), it can never go down too fast, thanks to Ray's estimate. Suppose that for some time $s\in[t,t+\t]$, $m^N_s$ is below $\g_s-7\ee$. Since $\t$ is an $\ee$-modulus of continuity for $\g$, we have that $\g_s-7\ee<\g_t-6\ee<r_5$. This means that on the event $\Xi_{r_4,r_5}^c$ with no downwards crossing of the interval $[r_5,r_4]$, a particle that lies below $\g_s-7\ee$ at time $s$ had to be below $r_4$ at time $t$, so that 
$$m^N_t\leq r_4<\g_t-4\ee.$$
Hence we have
$$\P(\exists s\in[t,t+\t],m^N_s<\g_s-7\ee)\leq\P(\Xi_{r_4,r_5})+\P_{\mu_0}(m^N_t<\g_t-4\ee).$$
Since the bound \eqref{eq:bound-probXi} holds as well for $\P(\Xi_{r_4,r_5})$, we deduce, using \eqref{eq:borne-t-fixe}, that 
 for every $\varkappa' > 0$, there exists a constant $C^{(2)}$ (depending on $\b,\varkappa',\ee,\g,s_0,t_0,T$) such that, for all $N \geq 1$, 
 \begin{equation}\label{eq:borne-t-petit-intervalle} \P(\exists s\in[t,t+\t],m^N_s<\g_s-7\ee)\leq C^{(2)} N^{-\varkappa'}.\end{equation} 
Taking the union bound over every possible $t = k \t$ in $[t_0-\t,T]$, where there are at most $\frac T\t$ of these intervals in $[t_0,T]$, we have  
  $$\P(\exists s\in[t_0,T],\ m^N_s<\g_s-7\ee)\leq\dfrac{T}{\t}  C^{(2)} N^{-\k}=2TC^{(2)} N^{\b-\varkappa'}.$$
  Given $\varkappa > 0$, we can set $\varkappa' = \b + \varkappa$, and deduce a bound, valid for all $N \geq 1$, of the form
\begin{equation}\label{eq:min-lowerbound} \P(\exists s\in[t_0,T],\ m^N_s<\g_s-7\ee)\leq C^{(3)} N^{-\varkappa},\end{equation}
 with a constant $C^{(3)}$ depending on $\b,\varkappa,\ee,\g,s_0,t_0,T$.

\noindent \textit{Step 4. Upper bound and conclusion. }

Now we prove that $\g$ is an approximate upper bound for the minimum for large $N$. We start by proving \eqref{eq:min-upperbound}, which bounds $\P(m^N_t>\g_t+3\ee)$ for $t=k\t$ in $[t_0,T]$.
Write $p:=\e^t U(\g_t+3\ee,t)$, consider $\b<\b'<1/2$, and let us now prove that, for large enough $N$ (depending on $\b,\b',\ee, t_0,T,\g$), we have \begin{equation}\label{eq:borne-p}p \geq 2 N^{-\b'}.\end{equation} 
For $i \in \Z$, let $r'_i=\min\{x_j \mbox{ such that } x_j>\g_t\}+i\ee$ and note that $r'_0,\ldots, r'_4$ are elements of the space grid subject to \eqref{eq:ray-grid} 
We then have $r'_2\leq \g_t+3\ee$ and thus $\e^t U(r'_2,t)\leq p$ since $r \mapsto U(r,t)$ is a non-decreasing function.  Since $\tau^\g$ follows the exponential distribution $\EE(1)$ under $\Q_{\mu_0}$, we have that 
$\Q_{\mu_0}(\tau^\g\leq t+\t|\tau^\g>t) = 1-\e^{-\t}$. On the other hand, 
$$\Q_{\mu_0}(\tau^\g\leq t+\t|\tau^\g>t)=\Q_{\mu_0}(X_t>r'_2,\tau^\g\leq t+\t | \tau^\g>t)+\Q_{\mu_0}(X_t\leq r'_2,\tau^\g\leq t+\t | \tau^\g>t).$$
With our assumptions, we have that $\g_s \leq r'_1$ for every $s \in [t,t+\t]$. Using the fact that that $\tau^{\g}$ is a stopping time, the Markov property, and stochastic monotonicity, we deduce that 
$\Q_{\mu_0}(X_t>r'_2,\tau^\g\leq t+\t | \tau^\g>t) \leq \Q_{r'_2}(\tau_{r'_1} \leq \t)$, and Ray's estimate \eqref{eq:ray-grid} ensures that this probability is bounded above by $C_\k\t^\k$ (for any $\k>0$). On the other hand, we have that $\Q_{\mu_0}(X_t\leq r'_2,\tau^\g\leq t+\t | \tau^\g>t) \leq \e^t U(r'_2,t) \leq p$. Putting together the previous results, we have proved that $1-\e^{-\t}\leq C_\k\t^\k+p$, so that $p \geq 1-\e^{-\t}-C_\k \t^\k$.
Remembering that $\t=2N^{-\b}$, and choosing $\k > 1$, \eqref{eq:borne-p} is proved.

Using the fact that  $p=\e^t U(\g_t+3\ee,t)$, we have that
\begin{equation}
\label{eq:min-upperbound}
    \P(m^N_t>\g_t+3\ee)=\P(F^N(\g_t+3\ee,t)=0)\leq\P\left(\n{F^N(\cdot,t)-\e^tU(\cdot,t)}\infty>\dfrac{p}{2}\right).
\end{equation}
 Thanks to \eqref{eq:min-upperbound},  \eqref{eq:borne-p} and \eqref{eq:hydro-nbmp-quant}, we have that 
\begin{equation}\label{eq:upperbound-instant}\P(m^N_t>\g_t+3\ee)\leq \P\left(\n{F^N(\cdot,t)-\e^tU(\cdot,t)}\infty>N^{-\b'}\right)\leq c_2\e^{-c_1N^{\frac{1-2\b'}3}}.\end{equation}

Then we argue in a similar way as in Step 3. Consider  $t=k\t$ in $[t_0,T-\t]$. Using the fact that $\g_s \geq \g_t - \ee \geq r'_{-2}$ for all $s \in [t,t+\t]$, we see that, on the event $\{ \exists s\in[t,t+\t],m^N_s>\g_s+7\ee \}$, the event $\{  m^N_{t+\t} \leq \g_{t+\t}+3\ee \}$ implies that at least one particle lies above $r'_5$ at a certain time $s \in [t,t+\t]$ then crosses $r'_4$ before time $t+\t$. Using Ray's estimate \eqref{eq:ray-cross-strong} and arguing exactly as in the proof of \eqref{eq:bound-probXi}, we deduce that the probability of the intersection of these two events is bounded above by $4 C_{\k} N^{1-\b\k}+\e^{-c_6N}$.

 In view of \eqref{eq:upperbound-instant}, we deduce that 
$$\P(\exists s\in[t,t+\t],m^N_s>\g_s+7\ee) \leq c_2\e^{-c_1N^{\frac{1-2\b'}3}} + 4 C_{\k} N^{1-\b\k}+\e^{-c_6N}.$$ Taking the union bound over every possible $t = k \t$ in $[t_0,T-\t]$ (whose number is bounded above e.g. by $T/\t$), we obtain an upper bound of the form 
\begin{equation}\label{eq:borne-sup-intervalle}\P(\exists s\in[t_0,T],m_s^N>\g_s+7\ee) \leq  C^{(4)} N^{-\k},\end{equation} 
valid for all $N \geq 1$ and $\varkappa > 0$.

Combining \eqref{eq:borne-sup-intervalle} and with \eqref{eq:min-lowerbound}, we deduce \eqref{eq:min}.

\noindent \textit{The case $t_0=0$. }

In the case $t_0=0$, the map $\g$ is continuous on the whole interval $[0,T]$, and we can take $s_0=0$ in the definition of the grid in Step 1. The proof of \eqref{eq:borne-t-fixe} for $t=k \theta$ with $k \geq 2$ in Step 2 requires no extra argument from the case $t_0>0$. On the other hand, when $k=0$, the proof is a consequence of the fact that, thanks to condition (i), we must have $\mu_0((-\infty,\g_0))=0$, so that $\P(m^N_0 \geq \g_0) = 1$. And when $k=1$, i.e. for $t=\theta$, this property also provides the control that is needed at time $t-\theta$ to prove \eqref{eq:borne-t-fixe} following the proof given in Step 2.  Step 3 is identical to the case $t_0>0$, and so is Step 4.

\noindent \textit{Almost sure convergence. }

Thanks to the Borel-Cantelli lemma, the almost sure convergence on $m^N_t$ to $\g_t$ in the supremum norm is a straightforward consequence of \eqref{eq:min}.

\end{demo}

\begin{rem}
As we mentioned earlier, the constants that appear in this proof are not universal: they depend strongly on the underlying Markov process $X$ through Ray's estimate and on the boundary $\g$ through the size of the grid and the number of terms used in the union bounds.
\end{rem}

\appendix

\section{Some results on the geometric distribution}
\label{sec:geom}

In this last section, we show some bounds on unlikely events regarding geometric distributions that will be used in several places in the proofs. 

\begin{lemma}
\label{lem:geom}
Let $G_1,\cdots,G_n$ be i.i.d. random variables with geometric distribution $\G(p)$, $p\in(0,1)$, and fix $q=1-p$ and $c=-\ln(q)>0$ . Then 
\begin{enumerate}[label=(\roman*)]
    \item for $K>0$, we have $\P(G_1>K)= \e^{-cK},$
    \item for $K>0$, we have $\E(G_1\1_{\{G_1>K\}})\leq (K+2)\dfrac{\e^{-cK}}{p}$,
    \item the Legendre transform of the geometric distribution $\G(p)$ is given, for $x\in(1,+\infty)$, by 
    \begin{equation}
    \label{eq:legendre}
        \Ll^*_p(x)=(x-1)\ln\left(\dfrac{x-1}{xq}\right)-\ln(px)\in[0,+\infty),
    \end{equation}
    and in particular, we have the following expansion for the Legendre transform $\Ll^*_\t$ when $p=\e^{-\t}$
    \begin{equation}
    \label{eq:expansion-leg}
        \Ll^*_{\e^{-\t}}\left(\frac1{1-\frac\t2}\right)=\dfrac{1-\ln(2)}{2}\t+o(\t),
    \end{equation}
    \item  for $x\in(1,\frac1p)$ and $n\geq 1$
    $$\P\left(\dfrac{1}{n}\dsum_{k=1}^nG_k<x\right)\leq\e^{-n\Ll^*_p(x)}.$$
\end{enumerate}
\end{lemma}

\begin{demo}
(i) is classical, and we can compute the exact value for (ii): we have 
$$\E(G_1\1_{\{G_1>K\}})=\dsum_{k=K+1}^{+\infty}kpq^{k-1}.$$
Set for $x<1$
$$g(x)=\dsum_{k=K+1}^{+\infty}pq^{k-1}x^k=pq^K\dfrac{x^{K+1}}{1-qx},$$
so that
$$\dsum_{k=K+1}^{+\infty}kpq^{k-1}x^{k-1}=g'(x)=pq^K\dfrac{(K+1)x^K(1-qx)+qx^{K+1}}{(1-qx)^2}.$$
Evaluating in $x=1$, we obtain 
$$\E(G_1\1_{\{G_1>K\}})=\dfrac{q^K((K+1)p+q)}{p}\leq (K+2)\dfrac{\e^{-cK}}{p}.$$
Let us carry out the computation of the Legendre transform of the geometric distribution: we have for $a<-\ln(q)$
$$\ln(\E(\e^{a G_1}))=\dsum_{k=1}^{+\infty} \e^{a k}q^{k-1}p=\dfrac{p\e^a}{1-q\e^a},$$
so we can compute for $x$ in $(1,+\infty)$ the derivative of $a\longmapsto a x-\ln(\E(\e^{a G_1}))$, which vanishes when $x-\frac{1}{1-q\e^a}$ is zero, \ie\ for $a^*=\ln\left(\frac{x-1}{xq}\right)$. In the end, we have as claimed
$$\Ll^*_p(x)=a^*x-\ln(\E(\e^{a^*G_1}))=(x-1)\ln\left(\dfrac{x-1}{xq}\right)-\ln(px).$$
The large deviations estimate (iv) then follows directly from Cramér's theorem (see e.g. \cite{DemZei}). 

There only remains to prove the expansion \eqref{eq:expansion-leg}, which we use in the proof of \thref{th:min}. Take $x=\frac{1}{1-\frac\t2}$ and $p=\e^{-\t}$ in \eqref{eq:legendre}. We have 
    $$\ln(x\e^{-\t})=\ln(x)-\t=-\dfrac\t2+o(\t).$$
    Next $\frac{x-1}{x}=\frac\t2$, so 
    $$\ln\left(\dfrac{x-1}{x(1-\e^{-\t})}\right)=\ln\left(\dfrac{\t}{2\t-\t^2+o(\t^2)}\right)=-\ln(2)+\dfrac\t2+o(\t),$$
    and then 
    $$(x-1)\ln\left(\dfrac{x-1}{x(1-\e^{-\t})}\right)=-\dfrac{\ln(2)}{2}\t+o(\t).$$
    In the end, we have indeed
    {\renewcommand\theequation{\ref{eq:expansion-leg}}
    \addtocounter{equation}{-1}
    \begin{equation}
        \Ll_{\e^{-\t}}^*\left(\frac1{1-\frac\t2}\right)=\dfrac{1-\ln(2)}{2}\t+o(\t).
    \end{equation}
    }
\end{demo}


\bibliographystyle{abbrv}
\bibliography{nbmp}

\begin{thebibliography}{10}

\bibitem{atar23}
R.~Atar.
\newblock Hydrodynamics of particle systems with selection via uniqueness for
  free boundary problems.
\newblock preprint, arXiv:2011.07535 [math.PR], 2023.

\bibitem{BerFre}
J.~B\'erard and B.~Fr\'enais.
\newblock Monotone coupling of {F}eller {M}arkov processes on the real line.
\newblock (forthcoming), 2023.

\bibitem{berard-gouere}
J.~B{\'e}rard and J.-B. Gou{\'e}r{\'e}.
\newblock {Brunet-Derrida behavior of branching-selection particle systems on
  the line}.
\newblock {\em {Communications in Mathematical Physics}}, 298(2):323--342,
  2010.

\bibitem{bees}
J.~Berestycki, E.~Brunet, J.~Nolen, and S.~Penington.
\newblock Brownian bees in the infinite swarm limit.
\newblock {\em Ann. Probab.}, 50(6):2133--2177, 2022.

\bibitem{fbp-fkpp}
J.~Berestycki, E.~Brunet, and S.~Penington.
\newblock Global existence for a free boundary problem of {F}isher-{KPP} type.
\newblock {\em Nonlinearity}, 32(10):3912--3939, 2019.

\bibitem{brunet-derrida}
E.~Brunet and B.~Derrida.
\newblock Shift in the velocity of a front due to a cutoff.
\newblock {\em Phys. Rev. E (3)}, 56(3):2597--2604, 1997.

\bibitem{chen22}
X.~Chen, J.~Chadam, and D.~Saunders.
\newblock Higher-order regularity of the free boundary in the inverse
  first-passage problem.
\newblock {\em SIAM J. Math. Anal.}, 54(4):4695--4720, 2022.

\bibitem{chen11}
X.~Chen, L.~Cheng, J.~Chadam, and D.~Saunders.
\newblock Existence and uniqueness of solutions to the inverse boundary
  crossing problem for diffusions.
\newblock {\em Ann. Appl. Probab.}, 21(5):1663--1693, 2011.

\bibitem{qsd}
P.~Collet, S.~Mart\'{\i}nez, and J.~San~Mart\'{\i}n.
\newblock {\em Quasi-stationary distributions}.
\newblock Probability and its Applications (New York). Springer, Heidelberg,
  2013.
\newblock Markov chains, diffusions and dynamical systems.

\bibitem{nbbm}
A.~De~Masi, P.~A. Ferrari, E.~Presutti, and N.~Soprano-Loto.
\newblock Hydrodynamics of the {N-BBM} process.
\newblock In G.~Giacomin, S.~Olla, E.~Saada, H.~Spohn, and G.~Stoltz, editors,
  {\em Stochastic Dynamics Out of Equilibrium}, pages 523--549, Cham, 2019.
  Springer International Publishing.

\bibitem{nonlocal-nbbm}
A.~De~Masi, P.~A. Ferrari, E.~Presutti, and N.~Soprano-Loto.
\newblock Non local branching {B}rownian motions with annihilation and free
  boundary problems.
\newblock {\em Electron. J. Probab.}, 24:Paper No. 63, 30, 2019.

\bibitem{DemZei}
A.~Dembo and O.~Zeitouni.
\newblock {\em Large deviations techniques and applications}, volume~38 of {\em
  Stochastic Modelling and Applied Probability}.
\newblock Springer-Verlag, Berlin, 2010.
\newblock Corrected reprint of the second (1998) edition.

\bibitem{devroye}
L.~Devroye, L.~Gy\"{o}rfi, and G.~Lugosi.
\newblock {\em A probabilistic theory of pattern recognition}, volume~31 of
  {\em Applications of Mathematics (New York)}.
\newblock Springer-Verlag, New York, 1996.

\bibitem{dudley}
R.~M. Dudley.
\newblock {Central Limit Theorems for Empirical Measures}.
\newblock {\em The Annals of Probability}, 6(6):899 -- 929, 1978.

\bibitem{durrett-remenik}
R.~Durrett and D.~Remenik.
\newblock Brunet-{D}errida particle systems, free boundary problems and
  {W}iener-{H}opf equations.
\newblock {\em Ann. Probab.}, 39(6):2043--2078, 2011.

\bibitem{qsd-levy}
P.~Groisman and M.~Jonckheere.
\newblock Front propagation and quasi-stationary distributions for
  one-dimensional {L}\'{e}vy processes.
\newblock {\em Electron. Commun. Probab.}, 23:Paper No. 93, 11, 2018.

\bibitem{groisman20}
P.~Groisman, M.~Jonckheere, and J.~Mart\'{\i}nez.
\newblock F-{KPP} scaling limit and selection principle for a
  {B}runet-{D}errida type particle system.
\newblock {\em ALEA Lat. Am. J. Probab. Math. Stat.}, 17(1):589--607, 2020.

\bibitem{groisman21}
P.~Groisman and N.~Soprano-Loto.
\newblock Rank dependent branching-selection particle systems.
\newblock {\em Electron. J. Probab.}, 26:Paper No. 158, 27, 2021.

\bibitem{harris-roberts}
S.~C. Harris and M.~I. Roberts.
\newblock The many-to-few lemma and multiple spines.
\newblock {\em Ann. Inst. Henri Poincar\'{e} Probab. Stat.}, 53(1):226--242,
  2017.

\bibitem{inw}
N.~Ikeda, M.~Nagasawa, and S.~Watanabe.
\newblock {Branching Markov processes I}.
\newblock {\em Journal of Mathematics of Kyoto University}, 8(2):233 -- 278,
  1968.

\bibitem{ito-mckean}
K.~It{\^o} and H.~P. McKean.
\newblock {\em Diffusion Processes and their Sample Paths}.
\newblock Classics in Mathematics. Springer Berlin Heidelberg, 2012.

\bibitem{kallenberg}
O.~Kallenberg.
\newblock {\em Foundations of Modern Probability}.
\newblock Springer, 2nd edition, 2002.

\bibitem{klump23}
A.~Klump and M.~Savov.
\newblock Conditions for existence and uniqueness of the inverse first-passage
  time problem applicable for {L}\'evy processes and diffusions.
\newblock preprint, arXiv:2305.10967 [math.PR], 2023.

\bibitem{kyprianou-palmowski}
A.~E. Kyprianou and Z.~Palmowski.
\newblock Quasi-stationary distributions for {L}\'{e}vy processes.
\newblock {\em Bernoulli}, 12(4):571--581, 2006.

\bibitem{lee}
J.~Lee.
\newblock A free boundary problem with non local interaction.
\newblock {\em Math. Phys. Anal. Geom.}, 21(3):Paper No. 24, 21, 2018.

\bibitem{maillard-nbbm}
P.~Maillard.
\newblock Speed and fluctuations of {$N$}-particle branching {B}rownian motion
  with spatial selection.
\newblock {\em Probab. Theory Related Fields}, 166(3-4):1061--1173, 2016.

\bibitem{qsd-bm}
S.~Mart\'{\i}nez and J.~San~Mart\'{\i}n.
\newblock Quasi-stationary distributions for a {B}rownian motion with drift and
  associated limit laws.
\newblock {\em J. Appl. Probab.}, 31(4):911--920, 1994.

\bibitem{pollard}
D.~Pollard.
\newblock {\em Convergence of {Stochastic Processes}}.
\newblock Springer, 1984.

\bibitem{revuz-yor}
D.~Revuz and M.~Yor.
\newblock {\em Continuous Martingales and Brownian Motion}.
\newblock Springer, 3rd edition, 1999.

\bibitem{yamato}
K.~Yamato.
\newblock Existence of quasi-stationary distributions for spectrally positive
  {L}\'{e}vy processes on the half-line.
\newblock {\em ALEA Lat. Am. J. Probab. Math. Stat.}, 20(1):629--643, 2023.

\end{thebibliography}

\end{document}